# Data-adaptive harmonic spectra and multilayer Stuart-Landau models

Mickaël D. Chekroun and Dmitri Kondrashov
*Department of Atmospheric and Oceanic Sciences, University of California, Los Angeles, CA 90095-1565, USA*

Harmonic decompositions of multivariate time series are considered for which we adopt an integral operator approach with periodic semigroup kernels. Spectral decomposition theorems are derived that cover the important cases of two-time statistics drawn from a mixing invariant measure.

The corresponding eigenvalues can be grouped per Fourier frequency, and are actually given, at each frequency, as the singular values of a cross-spectral matrix depending on the data. These eigenvalues obey furthermore a variational principle that allows us to define naturally a multidimensional power spectrum. The eigenmodes, as far as they are concerned, exhibit a data-adaptive character manifested in their phase which allows us in turn to define a multidimensional phase spectrum.

The resulting data-adaptive harmonic (DAH) modes allow for reducing the data-driven modeling effort to elemental models stacked per frequency, only coupled at different frequencies by the same noise realization. In particular, the DAH decomposition extracts time-dependent coefficients stacked by Fourier frequency which can be efficiently modeled—provided the decay of temporal correlations is sufficiently well-resolved—within a class of multilayer stochastic models (MSMs) tailored here on stochastic Stuart-Landau oscillators.

Applications to the Lorenz 96 model and to a stochastic heat equation driven by a space-time white noise, are considered. In both cases, the DAH decomposition allows for an extraction of spatio-temporal modes revealing key features of the dynamics in the embedded phase space. The multilayer Stuart-Landau models (MSLMs) are shown to successfully model the typical patterns of the corresponding time-evolving fields, as well as their statistics of occurrence.

## I. INTRODUCTION

The success of machine learning algorithms generally depends on data representation, and it is often hypothesized that this is because different representations can entangle and hide more or less the different explanatory factors of variations behind the data [1]. To be reasonably addressed, the problem of finding relevant data-adaptive basis functions for a successful inverse modeling, must be considered within a specific domain of applications [2]. This article is concerned with addressing this problem for time-evolving datasets produced by a dynamical system either deterministic or stochastic.

Among the numerous techniques useful for the data representation or decomposition of time-evolving datasets, we may distinguish several seminal techniques like those based on variance's decomposition such as principal component analysis (PCA) [3] and its probabilistic formulation [4] or nonlinear extensions [5–7]; techniques using eigenfunctions of Markov matrices reflecting the (local) geometry and density of the data [8, 9]; or more recent approaches exploiting the Koopman operator theory [10–12].

These methods often tied to dimensionality reduction, extract data-adaptive modes that come with reduction coordinates, which for datasets evolving in time, correspond to time series known as principal components (PCs) and the like. The effective derivation of reduced evolution equations to characterize the long-term dynamics of the underlying dynamical system[13] based on these coordinates, is then a central issue.

For some systems, the form of the master equations, prior experience, or some underlying physical intuition may help determine good reduced systems. In some special cases, their exact form can even be found rigorously by adopting, for instance, the Mori-Zwanzig (MZ) projection approach [14–19] or techniques rooted in the approximation theory of (stochastic) local invariant manifolds; see [20, 21] and references therein. In more complex cases where a rigorous derivation of the reduced dynamics is mathematically intractable, or when the master equations are even not known, a data-driven inverse modeling approach is often adopted.

Many methods have been proposed in the literature to address this task, such as those using: nonlinear autoregression moving average with exogenous inputs (NARMAX), artificial neural networks (ANNs), stochastic linear inverse models (LIMs), empirical model reduction (EMR), and many others; see [22–29] and references therein. In a data-driven context, the MZ formalism provides also a theoretical guidance for the determination of nonlinear ingredients and other dependences on the past; see [30]. Nevertheless, whether the approach retained uses an a priori set of predictor functions or allow for dictionaries of such functions [31], the modeler is often left, either with an explanatory deficit or with a selection problem of the appropriate (class of) predictors for a reliable emulation of the dynamics. In the latter case, a certain level of familiarity with the dataset is typically required to envision the right form of the predictors.

This article proposes to overcome such hurdles for inverse modeling, by introducing a data-adaptive spectral representation of the original data, aimed at providing—for a broad class of time-evolving datasets—reduction coordinates that can be efficiently modeled within a family of stochastic differential equations (SDEs) constituted from a fixed set of predictor functions; see Sectns. VII and VIII. Roughly speaking, the proposed modeling approach operates well for finite multivariate time series for which decays of correlations are sufficiently well-resolved for the amount of available data. It covers thus a broad class of datasets issued from dynamical systems that are mixing.

The approach adopted in this article for the data representation of time-evolving datasets relies on the spectral analysis of integral operators $\mathfrak{L}$ for which the kernels are built from periodic semigroups applied to the data's correlations. The resulting class of operators constitute one the two main contributions of this article, and a substantial portion of the latter is devoted to the analysis of their spectral properties; see Sectns. II, III, IV, and V. The other main contribution is concerned with the aforementioned class of SDEs, aimed at the modeling of the reduction coordinates obtained as projection of the original dataset onto the eigenfunctions of $\mathfrak{L}$. These SDEs fall into the class of networks of linearly coupled Stuart-Landau oscillators [32], that may be embedded within the class of multilayer stochastic models (MSMs) [30], allowing thus for the inclusion of memory effects in their formulation; see Sect. VII.

Another originality of the framework lies in its harmonic flavor. The framework allows indeed for extracting modes naturally ranked per Fourier frequency, with reduction coordinates that are, for mixing dynamical systems, typically narrowband and modulated in amplitude; see Sect. VII A. The corresponding reduced SDEs are also organized per Fourier frequency and the modeling of a specific frequency band is thus naturally made available and facilitated by the resulting framework.

The article is organized as follows. In Sect. II, we first recall basic elements of the spectral theory of periodic semigroups on a Banach space. In Sect. III, this spectral theory is applied to analyze the spectral properties of a class of integral operators $\mathscr{L}_\phi$ with periodic kernels built from periodic semigroups acting on a given (one-dimensional) data-function $\phi$; see Theorem III.1. The important case for applications of kernels built from left circular translation groups is analyzed in Sect. IV. Here, one recovers the Fourier spectrum (power and phase) from the spectral elements of $\mathscr{L}_\phi$ for which the eigenfunctions are phase-shifted sinusoids and their phase relates to the phase spectrum of $\phi$; see Theorem IV.1. This representation of the Fourier spectrum via the operator $\mathscr{L}_\phi$ allows us to propose in turn an innovative generalization for multidimensional signals in Sect. V. There we show, for multidimensional data, that the eigenvalues of an operator matrix generalization, $\mathfrak{L}$, can be grouped per Fourier frequency, and are actually given, at each frequency, as the singular values of a cross-spectral matrix depending on the data; see Theorem V.1. These eigenvalues obey furthermore a variational principle that allows us to define naturally a multidimensional power spectrum; see Remark V.1-(ii). The eigenmodes, as far as they are concerned, exhibit a data-adaptive character manifested in their phase which allows us in turn to define a multidimensional phase spectrum; see (81) below.

In Sect. VI we then consider the case for which $\phi$ is—or a collection of such $\phi$'s are—obtained from lagged correlations estimated from time series issued from a mixing dynamical system. After discussing a natural periodization to apply the theory (see Sect. VI B), we introduce in Sectns. VI C and VI D, the Hankel matrices and block-Hankel matrices resulting, respectively, from the discretization of $\mathscr{L}_\phi$ in the one-dimensional case, and of $\mathfrak{L}$ in the multidimensional case. First applications of the resulting data-adaptive harmonic (DAH) methodology, are then discussed for analytical examples of modulated and traveling waves; see Sect. VI E.

In Sect. VII, we introduce the class of MSMs tailored on stochastic Stuart-Landau oscillators, aimed at the modeling of the aforementioned reduction coordinates corresponding here to the projection of the dataset onto the DAH modes; see Sectns. VII A and VII B. These elemental models—the multilayer Stuart-Landau models (MSLMs)—are stacked per frequency, only coupled at different frequencies by the same noise realization; see Sect. VII C. In Sect. VIII, we finally show the flexibility of the DAH-MSLM modeling approach, by its ability to provide skilled stochastic inverse models for data issued either from the nonlinear chaotic Lorenz 96 model, or a stochastic heat equation driven by a space-time white noise. Concluding remarks discuss then in Sect. IX directions for future research.

## II. SPECTRAL THEORY OF PERIODIC SEMIGROUPS

Let us first recall the definition of a strongly continuous semigroup also known as $C_0$-semigroup; see [33, p. 36]. For an introduction to semigroup theory we refer to [34, 35]. More advanced treatments and applications of semigroup theory can be found in e.g. [33, 36–38].

**Definition II.1** *A family $(T(t))_{t\geq 0}$ of bounded linear operators on a Banach space $\mathscr{W}$ is called a strongly continuous (one-parameter) semigroup if*

$$\begin{aligned} &T(t+s) = T(t)T(s), \text{ for all } t,s \geq 0, \\ &T(0) = \text{Id}, \end{aligned} \tag{1}$$

*and for every $\varphi$ in $\mathscr{W}$, the orbit*

$$\begin{aligned} &\mathbb{R}^+ \to \mathscr{W} \\ &t \mapsto T(t)\varphi, \end{aligned} \tag{2}$$

*is continuous.*

The generator $A$ of a strongly continuous semigroup $(T(t))_{t\geq 0}$ is then defined as the operator

$$A : D(A) \subset \mathscr{W} \to \mathscr{W},$$

such that

$$A\varphi = \lim_{h\to 0^+} \frac{1}{h}\big(T(h)\varphi - \varphi\big), \tag{3}$$

for every $\varphi$ in the domain

$$D(A) := \{\varphi \in \mathscr{W} \mid \lim_{h\to 0^+} \frac{1}{h}\big(T(h)\varphi - \varphi\big) \text{ exists}\}. \tag{4}$$

Our approach rely on the following important characterization of spectral properties of periodic semigroups, e.g. [33,

Chap. IV, Sect. 2.25]. Recall that a strongly continuous semigroup $(T(t))_{t\geq 0}$ is periodic if there exists $t_0 > 0$ such that $T(t_0) = \mathrm{Id}$. The period $\tau$ of the semigroup is then

$$\tau = \inf\{t_0 > 0 : T(t_0) = \mathrm{Id}\}.$$

Periodic semigroups are always groups with inverses $T(t)^{-1} = T(n\tau - t)$, for $0 \leq t \leq n\tau$, and the spectrum of their generators lies always on the imaginary axis, more precisely we have

**Theorem II.1** *[33, Chap. IV, Lemma 2.25] Let $T(t)$ be a strongly continuous $\tau$-periodic semigroup on a Banach space $\mathscr{W}$, and A its generator, then*

$$\sigma(A) \subset \frac{2\pi i}{\tau}\cdot\mathbb{Z},$$

*and the resolvent $R(\mu,A)$ ($\mu \in \mathbb{C}$) of A is given by*

$$R(\mu,A) = (1-e^{-\mu\tau})^{-1}\int_0^\tau e^{-\mu s}T(s)\,\mathrm{d}s,\ \ \mu \notin \frac{2\pi i}{\tau}\cdot\mathbb{Z}$$

The above representation of the resolvent

$$R(\mu,A) := (A-\mu\mathrm{Id})^{-1},$$

shows that the resolvent of the generator of a $\tau$-periodic semigroup is a *meromorphic function* having only poles of maximal order, with poles $\mu_n$ given by

$$\mu_n = \frac{2i\pi n}{\tau}. \tag{5}$$

The residue of $R(\mu,A)$ at the pole $\mu_n$ is given by the bounded linear operator of $\mathscr{W}$ [33, p. 267],

$$P_n := \lim_{\mu\to\mu_n}(\mu-\mu_n)R(\mu,A) = \frac{1}{\tau}\int_0^\tau e^{-\mu_n s}T(s)\,\mathrm{d}s. \tag{6}$$

Hereafter, we denote by $\mathscr{L}(\mathscr{W})$ the space of bounded linear operators of the Banach space $\mathscr{W}$.

From [33, Chap. IV, §1.17], it follows furthermore that $P_n$ coincides with $\Pi_n$, the Riesz projector associated with the pole $\mu_n$, namely with

$$\Pi_n := -\frac{1}{2i\pi}\int_{\Gamma_c} R(\mu,A)\,\mathrm{d}\mu, \tag{7}$$

where $\Gamma_c$ denotes a rectifiable closed curve in the complex plane that separates $\mu_n$ from the rest of the poles. Furthermore the rank of $\Pi_n$, $\mathrm{rg}\,\Pi_n$, satisfies

$$\mathrm{rg}\,\Pi_n = \ker(A-\mu_n\mathrm{Id}),$$

i.e. the algebraic multiplicity is equal to the geometric multiplicity. In particular, this implies that the spectrum $\sigma(A)$ of $A$ consists of eigenvalues only (i.e. no continuous spectrum) and that

$$AP_n = \mu_n P_n. \tag{8}$$

Recall that for any strongly continuous $T(t)$ with generator $A$, for every $\lambda$ in $\mathbb{C}$, $t>0$, and $\varphi$ in $D(A)$, the following identity holds ([33, Chap. II, Lem. 1.9]):

$$e^{-\lambda t}T(t)\varphi - \varphi = \int_0^t e^{-\lambda s}T(s)(A-\lambda\mathrm{Id})\varphi\,\mathrm{d}s, \tag{9}$$

which here due to (8) leads to

$$T(t)P_n = e^{\mu_n t}P_n,\ \ t\geq 0. \tag{10}$$

Moreover,

$$\begin{aligned} P_mP_n\varphi &= \frac{1}{\tau}\int_0^\tau e^{-\mu_m s}T(s)\,\mathrm{d}s\cdot P_n\varphi \\ &= \frac{1}{\tau}\int_0^\tau e^{-\mu_m s}T(s)P_n\,\mathrm{d}s\cdot\varphi \\ &= \frac{1}{\tau}\int_0^\tau e^{(\mu_n-\mu_m)s}\,\mathrm{d}s\cdot P_n\varphi = 0,\ \text{for } n\neq m. \end{aligned} \tag{11}$$

With these properties at hand, it is then not difficult to infer the following useful decomposition result of any strongly continuous $\tau$-periodic semigroup and of its generator; see [33, Chap. IV, Thm. 2.27].

**Theorem II.2** *Let $T(t)$ be a strongly continuous $\tau$-periodic semigroup on a Banach space $\mathscr{W}$, and A its generator, then*

$$T(t)\varphi = \sum_{m\in\mathbb{Z}} e^{\mu_m t}P_m\varphi,\ \ \varphi\in D(A). \tag{12}$$

*and*

$$A\varphi = \sum_{m\in\mathbb{Z}}\mu_m P_m\varphi,\ \ \varphi\in D(A^2). \tag{13}$$

## III. SPECTRAL THEORY OF INTEGRAL OPERATORS WITH PERIODIC SEMIGROUP KERNELS

Let us introduce the one-dimensional torus $\mathbb{T} := \mathbb{R}/(\tau\mathbb{Z})$ which we define by identifying points in $\mathbb{R}$ that differ by $n\tau$ for some $n$ in $\mathbb{Z}$. We denote by $L^2(\mathbb{T})$ the space of complex-valued measurable functions on $\mathbb{R}$ that are $\tau$-periodic, and square-integrable with respect to an arbitrary reference interval of length $\tau$, $I = [\alpha, \alpha+\tau]$ for some $\alpha$ in $\mathbb{R}$. We endow $L^2(\mathbb{T})$ with the natural inner product

$$\langle f,g\rangle = \frac{1}{\tau}\int_{\mathbb{T}} f(r)\overline{g(r)}\,\mathrm{d}r,\ \ f,g\in L^2(\mathbb{T}). \tag{14}$$

Here the integral over $\mathbb{T}$ is the integral with respect to $r$ taken over any interval of length $\tau$. In what follows, we indeed often identify $\mathbb{T}$ with its interval representation $I$, making precise the distinction when necessary.

Given a $\tau$-periodic $C_0$-semigroup $T(t)$ acting on a function space $\mathscr{W}$ intersecting the set of $\tau$-periodic functions, and a $\tau$-periodic function $\phi$ in $\mathscr{W}$, our goal is to analyze the spectral property of the following operator:

$$\mathscr{L}_\phi[\Psi](r) := \frac{1}{\tau}\int_{\mathbb{T}}[T(r)\phi](s)\Psi(s)\,\mathrm{d}s,\ \ \Psi\in L^2(\mathbb{T}). \tag{15}$$

In practice the space $\mathscr{W}$ may differ from $L^2(\mathbb{T})$ (by requiring for instance other regularities on $\phi$), but to simplify the presentation we restrict ourselves to the case

$$\mathscr{W} = L^2(\mathbb{T}).$$

The following lemma results directly from the theory of integral operators

**Lemma III.1** *Assume $\phi$ lies in $L^2(\mathbb{T})$ and let $T(t)$ be a $\tau$-periodic $C_0$-semigroup on $L^2(\mathbb{T})$. If*

**(A1)** *$(r,s) \mapsto [T(r)\phi](s)$ belongs to $L^2(\mathbb{T}\times\mathbb{T})$,*

*then the operator $\mathscr{L}_\phi$ defined by (15), maps $L^2(\mathbb{T})$ into $L^2(\mathbb{T})$ and is* **compact.** *If furthermore*

**(A2)** *$\overline{[T(r)\phi](s)} = [T(s)\phi](r)$,*

*then the operator $\mathscr{L}_\phi$ is* **self-adjoint.**

**Proof.** The proof is standard and consists of noting that under the integrability condition **(A1)**, the operator $\mathscr{L}_\phi$ is a Hilbert-Schmidt operator on $L^2(\mathbb{T})$ and thus compact; see, e.g. [39, Thm. 6.12] or [40, Chap. 16,Theorem 3]. Finally, condition **(A2)** allows us to write

$$\begin{aligned}\langle \mathscr{L}_\phi[\Psi],\Psi\rangle &= \int_{\mathbb{T}}\int_{\mathbb{T}} [T(r)\phi](s)\Psi(s)\overline{\Psi(r)}\,\mathrm{d}s\,\mathrm{d}r\\ &= \overline{\int_{\mathbb{T}}\int_{\mathbb{T}} [T(s)\phi](r)\Psi(r)\overline{\Psi(s)}\,\mathrm{d}s\,\mathrm{d}r}\\ &= \overline{\langle \mathscr{L}_\phi[\Psi],\Psi\rangle},\end{aligned} \tag{16}$$

which shows that $\langle \mathscr{L}_\phi[\Psi],\Psi\rangle$ lies in $\mathbb{R}$ for any $\Psi$ in $L^2(\mathbb{T})$, and thus $\mathscr{L}_\phi$ is self-adjoint. ■

As a consequence of the spectral theorem for self-adjoint compact operators we have the following theorem, in which either $J=\{1,...,N\}$ with $N=\dim(\mathrm{rg}(\mathscr{L}_\phi))$ if $\mathscr{L}_\phi$ has finite rank, or $J=\mathbb{N}$ otherwise.

**Theorem III.1** *Let $\mathscr{L}_\phi$ defined in (15), with $T(t)$ and $\phi$ satisfying the assumptions of Lemma III.1 including* **(A1)** *and* **(A2)**. *Then 0 is an eigenvalue of $\mathscr{L}_\phi$ and there exist:*

*(a) countably many nonzero real numbers $\{\lambda_n\}_{n\in\mathbb{J}}$ either finitely many, or such that $\lambda_n \to 0$ if infinitely many, and*

*(b) an orthonormal basis $\{E_n\}_{n\in\mathbb{J}}$ of $cl(rg(\mathscr{L}_\phi))$, the topological closure of the range of $\mathscr{L}_\phi$ in $L^2(\mathbb{T})$,*

*for which*

$$\mathscr{L}_\phi[\Psi] = \sum_{n\in\mathbb{J}} \lambda_n \langle \Psi, E_n\rangle E_n. \tag{17}$$

*Each $\lambda_n$ is an eigenvalue and each $E_n$ is an eigenfunction.*

Due to the convergence to zero, we usually order the eigenvalues $\lambda_n$ in decreasing order according to absolute value:

$$|\lambda_1| \geq |\lambda_2| \geq \cdots. \tag{18}$$

The multiplicity of a given eigenvalue $\lambda$ is the number of times it is repeated in the list of eigenvalues given above, or, equivalently, it is the dimension of the $\lambda$-eigenspace $\ker(\mathscr{L}_\phi - \lambda\,\mathrm{Id})$. The multiplicity of any given eigenvalue is finite.

So far, we have applied classical results from the theory of integral operators to deal with the spectral properties of $\mathscr{L}_\phi$. We can derive furthermore a key identity satisfied by the eigenfunctions of $\mathscr{L}_\phi$ (see (21) below), by exploiting the spectral theory of periodic semigroups recalled in Sect. II and requiring more regularity on $\phi$.

For that purpose, we focus here on the spectral properties of $\mathscr{L}_\phi$ in $L^2_{\mathbb{R}}(\mathbb{T})$, i.e. the subspace of $L^2(\mathbb{T})$ constituted by real-valued functions. We have then the following result concerning the spectral elements of $\mathscr{L}_\phi$ in $L^2_{\mathbb{R}}(\mathbb{T})$, for which we make use of the following version of the Fourier transform $\widehat{\Psi}$ of a function $\Psi$ in $L^2(\mathbb{T})$:

$$\widehat{\Psi}(f) := \int_{\mathbb{T}} \Psi(s)e^{-2i\pi f s}\,\mathrm{d}s,\ \ f\in\mathbb{R}. \tag{19}$$

**Theorem III.2** *Assume $\phi$ lies in $D(A)$, the generator of a $\tau$-periodic $C_0$-semigroup $T(t)$ on $L^2(\mathbb{T})$. Assume that* **(A1)** *and* **(A2)** *hold. Assume*

**(A3)** *there exists an orthogonal basis $\{e_n\}$ of $L^2(\mathbb{T})$ for which $P_n$ defined in (6) is rank-one, i.e.*

$$P_n f = \langle f, e_n\rangle e_n. \tag{20}$$

*Let us denote by $E_n$ an eigenfunction of $\mathscr{L}_\phi$ in $L^2_{\mathbb{R}}(\mathbb{T})$ associated with an eigenvalue $\lambda_n$, then for any frequency $f_k = k/\tau$, we have*

$$\lambda_n \widehat{E}_n(f_k) = \langle \phi, e_k\rangle\langle e_k, E_n\rangle. \tag{21}$$

**Proof.** Let us consider an eigenfunction $E_n$ associated with an eigenvalue $\lambda_n$ of the operator $\mathscr{L}_\phi$ given in (15). Then the identity (12) of Theorem II.2 gives:

$$\begin{aligned}\tau\lambda_n e^{-\mu_k r}E_n(r) &= \int_{\mathbb{T}} e^{-\mu_k r}T(r)\phi(s)E_n(s)\,\mathrm{d}s,\\ &= \sum_{m\in\mathbb{Z}} e^{(\mu_m-\mu_k)r}\int_{\mathbb{T}} (P_m\phi)(s)E_n(s)\,\mathrm{d}s,\end{aligned} \tag{22}$$

which due to (5) leads, by integration in the $r$-variable, to

$$\tau\lambda_n \int_{\mathbb{T}} e^{-\mu_k r}E_n(r)\,\mathrm{d}r = \int_{\mathbb{T}} (P_k\phi)(s)E_n(s)\,\mathrm{d}s. \tag{23}$$

The latter identity gives (21) by applying (19) with $f=k/\tau$, combined with Assumption **(A3)** and the definition of the inner product $\langle\cdot,\cdot\rangle$ in (14), recalling that $E_n$ is real-valued. ■

**Remark III.1** *Under the conditions of Theorem III.1, one has furthermore*

$$\lambda_1 = \sup_{\Psi\in\mathscr{W},\Psi\neq 0} \frac{\langle \mathscr{L}_\phi\Psi,\Psi\rangle}{|\Psi|^2_{L^2}}, \tag{24}$$

*with $\mathscr{W}=L^2(\mathbb{T})$, e.g. [39, Prop. 6.9].*

*We can actually have a variational characterization of any eigenvalue of the self-adjoint compact operator $\mathscr{L}_\phi$ by relying on the Courant-Fischer min-max principle; see [39, Problem 37] and [41]. However since the operator $\mathscr{L}_\phi$ is not positive, the Courant-Fischer min-max principle needs to be amended. The positive eigenvalues of $\mathscr{L}_\phi$ listed in decreasing order and denoted by $0 \le \cdots \le \lambda_2^+ \le \lambda_1^+$, satisfy thus the relationships*

$$\lambda_{k+1}^+ = \min_{\substack{V \\ \mathrm{codim}(V)=k}} \max_{\substack{\Psi\in V \\ \Psi\neq 0}} \frac{\langle \mathscr{L}_\phi \Psi, \Psi\rangle}{|\Psi|_{L^2}^2}, \quad k\geq 0. \tag{25}$$

*That is, the* max *is taken over a linear subspace $V \subset \mathscr{W}$ of co-dimension k, and the* min *is taken over all such subspaces. Moreover, the minimum is attained on the subspace*

$$V = \mathrm{span}\{E_1^+, \cdots, E_k^+\}^\perp, \tag{26}$$

*where $E_j^+$ denotes the eigenfunction associated with a positive eigenvalue $\lambda_j^+$.*

*Similarly, for $0 \geq \cdots \geq \lambda_2^- \geq \lambda_1^-$ we have*

$$\lambda_{k+1}^- = \max_{\substack{V \\ \mathrm{codim}(V)=k}} \min_{\substack{\Psi\in V \\ \Psi\neq 0}} \frac{\langle \mathscr{L}_\phi \Psi, \Psi\rangle}{|\Psi|_{L^2}^2}. \tag{27}$$

*with the maximum attained for*

$$V = \mathrm{span}\{E_1^-, \cdots, E_k^-\}^\perp, \tag{28}$$

*where $E_j^-$ denotes the eigenfunction associated with a negative eigenvalue $\lambda_j^-$. Such min-max characterizations of eigenvalues of $\mathscr{L}_\phi$-like operators is made more explicit in Sect. V for the multidimensional case; see Remark V.1-(ii) below.*

## IV. THE CASE OF KERNELS FROM LEFT CIRCULAR TRANSLATION GROUPS

Let $g$ be a $\tau$-periodic continuous differentiable function $g$ such that $g'/g \in L^\infty$, and $\phi$ in $L^2(\mathbb{T})$. Let us consider

$$A\varphi = \frac{\mathrm{d}\varphi}{\mathrm{d}s} + \frac{g'}{g}\varphi, \tag{29}$$

whose domain is given by:

$$D(A) := \left\{ \varphi \in L^2(\mathbb{T}) \,:\, \varphi \in \mathrm{AC}(\mathbb{T}), \frac{\mathrm{d}\varphi}{\mathrm{d}s} + \frac{g'}{g}\varphi \in L^2(\mathbb{T}) \right\}, \tag{30}$$

where $\mathrm{AC}(\mathbb{T})$ denotes the collection of absolutely continuous functions on $\mathbb{T}$.

Then an adaptation of the arguments of [42, Lem. 3.5] shows that $A$ generates a weighted left translation semigroup $S(t)$ on $L^2(\mathbb{T})$. The corresponding operator $\mathscr{L}_\phi$ is then given by

$$\mathscr{L}_\phi[\Psi](r) = \frac{1}{\tau}\int_{\mathbb{T}} \frac{g(s+r)}{g(s)}\phi(s+r)\Psi(s)\,\mathrm{d}s.$$

By assumptions on $g$ and $\phi$, assumption (**A1**) is satisfied and $\mathscr{L}_\phi$ is compact; see Lemma III.1. If $g$ and $\phi$ are furthermore real-valued, then assumption (**A2**) is satisfied and thus $\mathscr{L}_\phi$ is self-adjoint as well; see again Lemma III.1.

An analysis of the resolvent of $A$ reveals that the corresponding operators $P_n$ defined in (6) are each rank-one, and that the $E_n$'s are real. Theorem III.2 can be thus applied here when $\phi$ lies in $D(A)$. Instead of providing a proof of this statement in this rather general setting, we focus on an important special case for applications, namely the case in which $g \equiv 1$, i.e. the case of the left circular translation group $T(t)$ defined by

$$T(t)\phi(s) = \phi(t+s), \quad \mod(\tau). \tag{31}$$

We have in this case

$$\mathscr{L}_\phi[\Psi](r) = \frac{1}{\tau}\int_{\mathbb{T}} \phi(s+r)\Psi(s)\,\mathrm{d}s, \tag{32}$$

and the corresponding Fredholm integral equation of the second kind[43] possesses explicit solutions such as described by the following theorem. These solutions provide actually a spectral representation of the Fourier transform in terms of both its power and phase spectra. This is the object of the next Theorem.

**Theorem IV.1** *Let $A$ be the generator of the left circular group $(T(t))_{t\in\mathbb{R}}$ on $L^2(\mathbb{T})$ defined by (31). Assume $\phi$ to be a real-valued function that lies in $D(A)$ given by (30) with $g \equiv 1$. Then the eigenvalues of the operator $\mathscr{L}_\phi$ are real-valued, form a discrete set $\{\lambda_k\}$, where $k$ in $\mathbb{Z}$ characterizes the Fourier frequency $2i\pi k/\tau$. More precisely, to each frequency $f_k = \frac{k}{\tau}$, corresponds an eigenvalue $\lambda_k$ such that*

$$|\lambda_k| = |\widehat{\phi}(f_k)|. \tag{33}$$

*Furthermore the corresponding eigenfunctions in $L^2_{\mathbb{R}}(\mathbb{T})$ are given by*

$$E_k(t) = \sqrt{2}\cos(2\pi f_k t + \theta_k), \tag{34}$$

*with*

$$\begin{cases} \theta_k = -\frac{1}{2}\arg\widehat{\phi}(f_k), \quad \mod(2\pi) \;\; \text{if } \lambda_k \geq 0, \\ \theta_k = -\frac{1}{2}\arg\widehat{\phi}(f_k) + \frac{\pi}{2}, \quad \mod(2\pi) \;\; \text{if } \lambda_k < 0, \end{cases} \tag{35}$$

*where* $\arg(z)$ *denotes the principal value[44] of the argument of the complex number z.*

**Remark IV.1** *The eigenfunctions are thus phase-shifted sinusoids, for which the phase $\theta_k$ relates to the phase spectrum of the data $\phi$. The eigenvalues provide the power spectrum. So far, for a one-dimensional signal $\phi$, one can recover the Fourier spectrum (power and phase) from the spectral elements of $\mathscr{L}_\phi$. The operator representation of the Fourier spectrum via the operator $\mathscr{L}_\phi$ allows us to propose an innovative generalization for multidimensional signals (see Sect. V) turns out to be particularly useful for the inverse stochastic-dynamic modeling of spatio-temporal scalar fields; see Sectns. VII and VIII.*

**Proof.** As explained above, assumptions **(A1)** and **(A2)** are satisfied. Let us now check assumption **(A3)** of Theorem III.2. For that let us determine the resolvent of the generator $A$ of the left circular group $(T(t))_{t\in\mathbb{R}}$ defined by (31). Note that the operator $A$ is the differentiation operator defined in (29) with $g \equiv 1$.

The resolvent $R(\lambda,A)$ of $A$ exists for all $\lambda \notin \{0, \pm i\omega_1, \pm i\omega_2, \cdots\}$, with $\omega_k = 2\pi k/\tau$, and is obtained by solving the differential equation

$$\lambda u - \dot{u} = \psi, \ \ \psi \in L^2(\mathbb{T}). \tag{36}$$

The variation-of-constant formula gives

$$u(s) = e^{\lambda(s-t)}u(t) - \int_t^s e^{\lambda(s-r)}\psi(r)\,\mathrm{d}r. \tag{37}$$

Taking $s = t+\tau$, and requiring $u(t+\tau) = u(t)$, we get

$$(1-e^{\lambda\tau})u(t) = -\int_t^{t+\tau} e^{\lambda(t+\tau-r)} f(r)\,\mathrm{d}r. \tag{38}$$

The factor on the left hand side is invertible if and only if $\lambda$ does not coincide with $2\pi i n/\tau$. In such a case,

$$\begin{aligned} u(t) &= (1-e^{-\lambda\tau})^{-1}\int_t^{t+\tau} e^{\lambda(t-r)}\psi(r)\,\mathrm{d}r \\ &= (1-e^{-\lambda\tau})^{-1}\int_0^{\tau} e^{-\lambda\tau}\psi(t+r)\,\mathrm{d}r. \end{aligned} \tag{39}$$

The right-hand side of this formula maps $L^2(\mathbb{T})$ into the Sobolev space $W^{1,2}(\mathbb{T})$.

Now let $p$ be arbitrary in $\mathbb{Z}$ and let us take $\psi(r) := e^{i\omega_p r}$. Then

$$\begin{aligned} u(t) &= e^{i\omega_p t}(1-e^{-\lambda\tau})^{-1}\int_0^{\tau} e^{(i\omega_p-\lambda)r}\,\mathrm{d}r \\ &= \frac{(e^{(i\omega_p-\lambda)\tau}-1)}{1-e^{-\lambda\tau}}\frac{e^{i\omega_p t}}{i\omega_p-\lambda} = -\frac{e^{i\omega_p t}}{i\omega_p-\lambda}, \end{aligned} \tag{40}$$

since $\omega_p\tau \in 2\pi\mathbb{Z}$, by definition. Thus by expanding now any arbitrary $\Psi$ in $L^2(\mathbb{T})$ into its Fourier expansion, we get

$$R(\lambda,A)\psi = \sum_{p=-\infty}^{\infty}\frac{\widehat{\psi}(f_p)}{i\omega_p-\lambda}e^{i\omega_p\theta}, \tag{41}$$

with $f_p = p/\tau$.

Therefore the residue $P_k$ at the pole $\mu_k = i\omega_k$ (see (6)) is given, for any $k$ in $\mathbb{Z}$, by

$$P_k\psi(t) = \widehat{\psi}(f_k)e^{i\omega_k t}, \tag{42}$$

and the assumption **(A3)** of Theorem III.2 is satisfied with $e_k(t) = e^{i\omega_k t}$.

The identity (21) becomes in this case

$$\lambda_k\widehat{E}_k(f_k) = \widehat{\phi}(f_k)\overline{\widehat{E}_k(f_k)}. \tag{43}$$

Taking the modulus on both sides of (43) gives

$$|\lambda_k| = |\widehat{\phi}(f_k)|. \tag{44}$$

Let us seek solutions to (43) under the form

$$\begin{aligned} E_k(t) &= \cos\Big(2\pi f_k t + \theta_k\Big) \\ &= \cos\Big(2\pi f_k\Big(t - \frac{\theta_k}{2\pi f_k}\Big)\Big). \end{aligned} \tag{45}$$

Now by using the Fourier time-shift property, i.e. that the Fourier transform of $x(t-\theta_k)$ is given by $\widehat{x}(f)e^{2\pi i f\theta}$, we obtain from the Fourier transform of $\cos(2\pi f_k t)$ that

$$\begin{aligned} \widehat{E_k}(f) = \frac{1}{2}\delta(f-f_k)\exp\Big(-\frac{if\theta_k}{f_k}\Big) \\ + \frac{1}{2}\delta(f+f_k)\exp\Big(-\frac{if\theta_k}{f_k}\Big). \end{aligned} \tag{46}$$

and

$$\widehat{E_k}(f_k) = \frac{1}{2}\exp\big(-i\theta_k\big). \tag{47}$$

With this expression substituted into (43), one finds

$$\lambda_k = \widehat{\phi}(f_k)\exp\Big(2i\theta_k\Big), \tag{48}$$

and since $\lambda_k$ is real, necessarily:

$$\begin{cases} \theta_k = -\frac{1}{2}\arg\widehat{\phi}(f_k), & \mod\big(2\pi\big) \ \text{ if } \lambda_k \geq 0, \\ \theta_k = -\frac{1}{2}\arg\widehat{\phi}(f_k) + \frac{\pi}{2}, & \mod\big(2\pi\big) \ \text{ if } \lambda_k < 0, \end{cases} \tag{49}$$

Finally, the factor $\sqrt{2}$ comes as a normalization constant of $E_k$, for the norm subordinated to the inner product (14). ∎

**Remark IV.2**

*(i) The negative part of the spectrum of $\mathscr{L}_\phi$ comes with eigenmodes of sine-type, namely if $\lambda_k < 0$ then*

$$\begin{aligned} E_k(t) &= \cos\big(\omega_k t - \frac{1}{2}\arg\widehat{\phi}(f_k) + \frac{\pi}{2}\big) \\ &= \sin\big(\omega_k t - \frac{1}{2}\arg\widehat{\phi}(f_k)\big). \end{aligned} \tag{50}$$

*(ii) Sometimes instead of (35), we will make use of the following equivalent expression of the phase spectrum in terms of the eigenmodes of $\mathscr{L}_\phi$, namely*

$$\begin{aligned} \arg\widehat{\phi}(f_k) = \arg\Big(\lambda_k\widehat{E}_k(f_k)\Big) \\ - \arg\Big(\overline{\widehat{E}_k(f_k)}\Big), \quad \mod(2\pi). \end{aligned} \tag{51}$$

*The latter identity is a direct consequence of Eq. (43).*

*(iii) If $\tau$ is a positive integer $N$, then the spectrum of $\mathscr{L}_\phi$ is finite and given by*

$$\sigma(\mathscr{L}_\phi) = \Big\{\big|\widehat{\phi}(f_k)\big| \ : \ k \in \{0,\cdots,N-1\}\Big\}. \tag{52}$$

## V. MULTIDIMENSIONAL GENERALIZATION

In this section we consider a block operator matrix generalization of what precedes. We assume that we are given a collection $\{\phi^{i,j}(t),\ 1 \leq i,j \leq d\}$ of $\tau$-periodic functions in $L^2(\mathbb{T})$. Examples of such functions for the spectral analysis of e.g. spatio-temporal fields is discussed in Sect. VI D below.

For the moment, let us consider the following operator matrix $\mathfrak{L}$ defined for each $u$ in $\left(L^2(\mathbb{T})\right)^d$ by:

$$\begin{aligned} \mathfrak{L}u &:= \left(\mathscr{M}_1\Psi_1, \cdots, \mathscr{M}_d\Psi_d\right), \\ u &= (\Psi_1, \cdots, \Psi_d) \in \left(L^2(\mathbb{T})\right)^d, \end{aligned} \tag{53}$$

with

$$\mathscr{M}_p\Psi_p := \sum_{q=1}^{p-1} \mathscr{L}_{\phi^{q,p}}(\Psi_q) + \sum_{q=p}^{d} \mathscr{L}_{\phi^{p,q}}(\Psi_q),\ \ 1 \leq p \leq d, \tag{54}$$

where each $\mathscr{L}_{\phi^{i,j}}$ is given by (32). We have then the following Lemma.

**Lemma V.1** *The operator $\mathfrak{L}$ is self-adjoint and compact on $\mathscr{X} := \left(L^2(\mathbb{T})\right)^d$ endowed with the inner product given for any $u = (\Psi_1, \cdots, \Psi_d)$ and $v = (v_1, \cdots, v_d)$ in $\mathscr{X}$ by*

$$\langle u, v\rangle_{\mathscr{X}} := \sum_{n=1}^{d} \langle \Psi_n, v_n\rangle = \tau^{-1} \sum_{n=1}^{d} \int_{\mathbb{T}} \Psi_n(r)\overline{v_n(r)}\,\mathrm{d}r. \tag{55}$$

**Proof.** Since for each $p,q$ in $\{1, \cdots, d\}$, the function $\phi^{p,q}$ lies in $L^2(\mathbb{T})$ by assumption, each operator $\mathscr{L}_{\phi^{p,q}}$ given by (32) is compact (see Sect. IV), and thus each operator $\mathscr{M}_p$ is a (real-valued) compact operator on $L^2(\mathbb{T})$, as the finite sum of compact operators. The operator $\mathfrak{L}$ is then compact on $\mathscr{X}$ as the finite cartesian product of compact operators, for $\mathscr{X}$ endowed with the product topology induced by the inner product (55).

The self-adjoint property of $\mathfrak{L}$ results from its definition. Indeed,

$$\langle \mathfrak{L}u, u\rangle_{\mathscr{X}} = \sum_{p=1}^{d} \Big\langle \sum_{q=1}^{p-1} \mathscr{L}_{\phi^{q,p}}(\Psi_q) + \sum_{q=p}^{d} \mathscr{L}_{\phi^{p,q}}(\Psi_q), \Psi_p \Big\rangle, \tag{56}$$

which, after similar manipulations as in (16), leads to,

$$\begin{aligned} \langle \mathfrak{L}u, u\rangle_{\mathscr{X}} &= \sum_{p=1}^{d} \Big( \sum_{q=1}^{p-1} \overline{\langle \mathscr{L}_{\phi^{q,p}}(\Psi_p), \Psi_q\rangle} + \sum_{q=p}^{d} \overline{\langle \mathscr{L}_{\phi^{p,q}}(\Psi_p), \Psi_q\rangle} \Big) \\ &= \sum_{q=1}^{d} \Big( \sum_{p=1}^{q-1} \overline{\langle \mathscr{L}_{\phi^{p,q}}(\Psi_p), \Psi_q\rangle} + \sum_{p=q}^{d} \overline{\langle \mathscr{L}_{\phi^{q,p}}(\Psi_p), \Psi_q\rangle} \Big) \\ &= \overline{\langle \mathfrak{L}u, u\rangle_{\mathscr{X}}}. \end{aligned} \tag{57}$$

■

We are now in position to formulate the main result that extends Theorem IV.1 to the multidimensional case. For that purpose, we introduce for each frequency $f_k = k/\tau$, the $d \times d$ **symmetrized cross-spectral matrix** $\Phi(f_k)$ whose entries are given by

$$\Phi^k_{p,q} = \begin{cases} \widehat{\phi^{p,q}}(f_k) \text{ if } q \geq p, \\ \widehat{\phi^{q,p}}(f_k) \text{ if } q < p. \end{cases} \tag{58}$$

The next theorem shows that the spectrum of the operator $\mathfrak{L}$ relates naturally to the symmetrized cross-spectral matrix $\Phi(f)$. At every frequency line, the singular values of $\Phi(f)$ gives actually $d$ pairs of eigenvalues of $\mathfrak{L}$, and as $f$ is varied one recovers the full spectrum of $\mathfrak{L}$[45].

**Theorem V.1** *Let $A$ be the generator of the left circular group $(T(t))_{t\in\mathbb{R}}$ on $L^2(\mathbb{T})$ defined by (31). Assume that each $\phi^{p,q}$ is a real-valued function lying in $D(A)$ given by (30) with $g \equiv 1$. Then the eigenvalues of the operator $\mathfrak{L}$ defined in (53), are real-valued, form a discrete set $\{\lambda_n\}$ that possesses the following characterization.*

*For each frequency $f_k \neq 0$, one can extract $2d$ eigenvalues (counting multiplicities) from the set $\{\lambda_n\}$, $d$ positive and $d$ negative ones that satisfy the following property.*

- **(S)** *For each pair of negative-positive eigenvalues denoted by $(\lambda^p_-(f_k), \lambda^p_+(f_k))$, there exists a* **singular value** *$\sigma_p(f_k)$* **of** *$\Phi(f_k)$ such that*

$$\lambda^p_+(f_k) = -\lambda^p_-(f_k) = \sigma_p(f_k),\ \ 1 \leq p \leq d. \tag{59}$$

*Furthermore the eigenfunctions $W_k$ of $\mathfrak{L}$ in $\left(L^2_{\mathbb{R}}(\mathbb{T})\right)^d$ possess the following representation*

$$W_k = (E^k_1, \cdots, E^k_d)^{\mathrm{tr}}, \tag{60}$$

*with*

$$\begin{aligned} E^k_p(t) = B^k_p \cos(2\pi f_k t + \theta^k_p), \quad B^k_p,\ \theta^k_p \in \mathbb{R}, \\ \text{with } f_k = \frac{k}{\tau},\ k \in \mathbb{Z},\ p \in \{1, \cdots, d\}. \end{aligned} \tag{61}$$

*Finally, introducing the notations,*

$$\widehat{\phi^{a,b}}(f_k) = X^k_{ab} e^{i\alpha^k_{ab}},\ \text{with } X^k_{ab} = |\widehat{\phi^{a,b}}(f_k)|, \tag{62}$$

*the amplitudes $B^k_p$ and angles $\theta^k_p$ satisfy for each frequency $f_k$ and each $p$ in $\{1, \cdots, d\}$, the relation*

$$\begin{aligned} B^k_p \lambda^p_+(f_k) = \sum_{q=1}^{p-1} B^k_q X^k_{qp} e^{i(\theta^k_p + \theta^k_q + \alpha^k_{qp})} + \\ \sum_{q=p}^{d} B^k_q X^k_{pq} e^{i(\theta^k_p + \theta^k_q + \alpha^k_{pq})}. \end{aligned} \tag{63}$$

**Proof. Step 1.** We first start with the derivation of an analogue of Eq. (43) in the multdimensional case considered here. The resulting extension of (43) is key to the derivation of the main properties of the spectral elements of the block operator matrix $\mathfrak{L}$ defined by (53)-(54).

For that purpose, let us consider $W_k$ in $\left(L^2_{\mathbb{R}}(\mathbb{T})\right)^d$ under the form given in (60) with the $E^k_p$ to be determined. Due to the operator matrix form of $\mathfrak{L}$ given by (53)-(54), the eigenvalue problem $\mathfrak{L}W_n = \lambda W_n$ can be equivalently rewritten as

$$\sum_{q=1}^{p-1} \mathscr{L}_{\phi^{q,p}}(E^k_q) + \sum_{q=p}^{d} \mathscr{L}_{\phi^{p,q}}(E^k_q) = \lambda_k E^k_p,\ \ p \in \{1, \cdots, d\}. \tag{64}$$

Let $p$ be fixed in $\{1,\cdots,d\}$ and $k$ be in $\mathbb{Z}$. By reproducing the arguments of the proof of Theorem III.2, we have:

$$\begin{aligned}\lambda \int_{\mathbb{T}} e^{-\mu_k r} E_p^n(r)\,\mathrm{d}r = &\sum_{q=1}^{p-1}\int_{\mathbb{T}} (P_k\phi^{q,p})(s)E_q^n(s)\,\mathrm{d}s \\ &+ \sum_{q=p}^{d}\int_{\mathbb{T}} (P_k\phi^{p,q})(s)E_q^n(s)\,\mathrm{d}s.\end{aligned} \tag{65}$$

Here again $P_k$ denotes the residue of the resolvent $R(\xi,A)$ of $A$ at the pole $\mu_k = 2\pi ik/\tau$, and is given thus by

$$P_k = \lim_{\xi\to\mu_k}(\xi-\mu_k)R(\xi,A).$$

Now let us recall the expression of the resolvent of $R(\lambda,A)$ derived in the proof of Theorem IV.1. Namely, $R(\xi,A)$ exists for all $\xi \notin \{0,\pm i\omega_1,\pm i\omega_2,\cdots\}$, with $\omega_k = 2\pi k/\tau$, and is given for $\psi$ in $L^2(\mathbb{T})$ by

$$R(\xi,A)\psi = \sum_{p=-\infty}^{\infty} \frac{\widehat{\psi}(f_p)}{i\omega_p - \xi} e^{i\omega_p \theta}. \tag{66}$$

Therefore $P_k$ is given by

$$P_k\psi(t) = \widehat{\psi}(f_k)e^{i\omega_k t}, \tag{67}$$

and thus from (65) we obtain, since $E_q^n(t)$ is a real-valued function (since $W_n$ lies in $\left(L^2_{\mathbb{R}}(\mathbb{T})\right)^d$),

$$\lambda \widehat{E_p^n}(f_k) = \sum_{q=1}^{p-1} \widehat{\phi^{q,p}}(f_k)\overline{\widehat{E_q^n}(f_k)} + \sum_{q=p}^{d} \widehat{\phi^{p,q}}(f_k)\overline{\widehat{E_q^n}(f_k)}, \tag{68}$$

i.e. the aforementioned extension of (43) to the multidimensional case.

Seeking solutions to Eq. (68) of the form (61), and recalling that

$$\begin{aligned}\widehat{E_p^n}(f) = &\frac{1}{2}\delta(f-f_n)\exp\left(-\frac{if\theta_p^n}{f_n}\right) \\ &+ \frac{1}{2}\delta(f+f_n)\exp\left(-\frac{if\theta_p^n}{f_n}\right).\end{aligned} \tag{69}$$

(see also (46)), we conclude that

$$(\widehat{E_p^n}(f_k)\neq 0 \text{ and } \widehat{E_q^n}(f_k)\neq 0) \iff (k=n). \tag{70}$$

We consider thus $k=n$, hereafter. Replacing $\lambda_+^p$ by $\lambda$ in Eq. (63) for a moment, it is easy to see that Eq. (63) is a direct consequence of Eq. (68). The former is indeed obtained by taking the real and imaginary parts of Eq. (68) (in which $E_p^k$ is given by (61)) and forming the complex exponentials by summation. To complete the proof, we are thus left with the proof of property **(S)** which we deal with next.

**Step 2.** Note first that the symmetrized cross-spectral matrix $\Phi(f_k)$ whose entries are defined in (58), is complex and symmetric. The Autonne-Takagi factorization (see [46, Chap. VIII, Théorème 92] and [47]) applied to $\Phi(f_k)$ ensures thus the existence of a unitary matrix $U_k$ such that

$$\Phi(f_k) = U_k\Sigma_k U_k^T. \tag{71}$$

where $\Sigma_k = \operatorname{diag}\{\sigma_1(f_k),\cdots,\sigma_d(f_k)\}$ is a diagonal matrix with real and nonnegative entries. These entries are the singular values of $\Phi(f)$ at the frequency $f=f_k$.

In particular, denoting by $U^*$ the conjugate transpose of $U$, we have

$$\begin{aligned}\Phi(f_k)\overline{\Phi(f_k)} &= U_k\Sigma U_k^T\overline{U}\Sigma_k U_k^* \\ &= U_k\Sigma^2 U_k^*,\end{aligned} \tag{72}$$

and therefore the matrix $\Phi(f_k)\overline{\Phi(f_k)}$ is positive semi-definite, and its eigenvalues are given by the $\sigma_j(f_k)^2$, $1\le j\le d$.

On the other hand, Eq. (68) may be rewritten using (60) as

$$\lambda\widehat{W_k}(f_k) = \Phi(f_k)\overline{\widehat{W_k}(f_k)} \tag{73}$$

which leads to

$$\begin{aligned}\Phi(f_k)\overline{\Phi(f_k)}\widehat{W_k}(f_k) &= \Phi(f_k)\overline{(\Phi(f_k)\overline{\widehat{W_k}(f_k)})} \\ &= \lambda\Phi(f_k)\overline{\widehat{W_k}(f_k)} \\ &= \lambda^2\widehat{W_k}(f_k),\end{aligned} \tag{74}$$

and thus $\lambda^2$ is an eigenvalue of $\Phi(f_k)\overline{\Phi(f_k)}$.

To summarize, $\lambda$ is (up to a sign) the square root of an eigenvalue of $\Phi(f_k)\overline{\Phi(f_k)}$ and is thus also a singular value $\sigma_j(f_k)$ of $\Phi(f_k)$, up to a sign. Now since $k=n$ in Eq. (68) (see step 1), there is for each frequency $f_k$, $d$ relations (68), one for each $p$. Considering the positive and negative eigenvalues that coexist at each frequency $f_k\neq 0$, we may infer thus (59). ■

**Remark V.1**

*(i) It should be noted that an SVD of $\Phi(f_k)$ is not automatically its Autonne-Takagi factorization. The proof of the Autonne-Takagi factorization provides the details of how an SVD of a complex symmetric matrix can be turned into a symmetric SVD; see the proof of [48, Theorem 3.1].*

*(ii) It is also worth mentioning that due to a general result dealing with the characterization of singular values of complex symmetric matrices [49], the singular values of $\Phi(f_k)$ possess a variational characterization in terms of the quadratic form*

$$\mathscr{Q}(x) := \operatorname{Re}(x^T\Phi(f_k)x),\ \ x\in\mathbb{C}^d. \tag{75}$$

*In particular by application of [49, Corollary 5] and Theorem V.1, we obtain for $0\le p<d$ the following variational characterization by maximizing over real subspaces, namely*

$$\max_{\substack{V\\ \dim_{\mathbb{R}}(V)=k+1}}\ \min_{\substack{x\in V\\ \|x\|=1}} \operatorname{Re}(x^T\Phi(f_k)x) = \lambda_+^{1+p}(f_k), \tag{76}$$

*and*

$$\max_{\substack{V\\ \operatorname{codim}_{\mathbb{R}}(V)=k}}\ \min_{\substack{x\in V\\ \|x\|=1}} \operatorname{Re}(x^T\Phi(f_k)x) = \lambda_-^{1+p}(f_k). \tag{77}$$

Theorem V.1 and Remark V.1-(ii) ensure that the **eigenvalues** of the operator $\mathfrak{L}$ (at the frequency $f = f_k$) relate naturally to an **energy via the quadratic form** $\mathscr{Q}$ given in (75). Contrarily to the one-dimensional case (see (35)) and as (63) shows, the $\theta_p^k$'s do not relate trivially to the multidimensional phase information contained in symmetrized cross-spectral matrix $\Phi(f_k)$.

Nevertheless, the relation (68) allows us to extract a useful phase information from $\Phi(f_k)$ by noting that (provided that $B_p^k \neq 0$)

$$\lambda_+^p(f_k)\widehat{E_p^k}(f_k)(\overline{\widehat{E_p^k}(f_k)})^{-1} = \widehat{\phi^{p,p}}(f_k) + R_p(f_k), \tag{78}$$

with

$$R_p(f_k) = \sum_{q=1}^{p-1} \widehat{\phi^{q,p}}(f_k)\gamma_{pq}^k + \sum_{q=p+1}^{d} \widehat{\phi^{p,q}}(f_k)\gamma_{pq}^k, \tag{79}$$

and $\gamma_{pq}^k := \overline{\widehat{E_q^k}(f_k)(\widehat{E_p^k}(f_k))^{-1}}$. Obviously a similar relation holds for $\lambda_-^p(f_k)$ (involving the appropriate $E_p^k$).

By taking the argument on both sides of (78), we obtain (up to a multiple of $2\pi$)

$$\begin{aligned}\arg(\lambda_+^p(f_k)\widehat{E_p^k}(f_k)) \\ -\arg(\overline{\widehat{E_p^k}(f_k)}) = \arg\left(\widehat{\phi^{p,p}}(f_k) + R_p(f_k)\right).\end{aligned} \tag{80}$$

We denote by $\eta_+^p(f_k)$ the LHS of (80).

On the other hand (69) gives $\widehat{E_p^k}(k) = \exp\left(-i\theta_p^k\right)/2$, which leads to

$$2\theta_p^k = \eta_+^p(f_k) = \arg\left(\widehat{\phi^{p,p}}(f_k) + R_p(f_k)\right). \tag{81}$$

We understand thus that $\eta_+^p(f_k)$ gives (twice) the phase shift, $\theta_p^k$, contained in the eigenfunctions (see (61)), and provides also a measure about the perturbation brought by the "cross-functions" $\phi^{p,q}$ ($p \neq q$)—as encapsulated into $R_p(f_k)$—to the phase spectrum associated with $\phi^{p,p}$ when considered as isolated (i.e. $\arg(\widehat{\phi^{p,p}}(f_k))$).

This perturbation is dependent on the eigenfunctions of the operator $\mathfrak{L}$ as the expression of $R_p(f_k)$ given in (79), shows. When $R_p(f_k) = 0$, on recovers the phase spectrum of $\phi^{p,p}$ as in the unidimensional case; see Remark VI.2-(iii). In what follows, we will look at the quantities, $\lambda_+^p(f_k)$ and $\eta_+^p(f_k)$ (resp. $\eta_+^p(f_k)$) as $f_k$ is varied; the former providing, in the multidimensional case, a notion of power spectrum (as associated with the form $\mathscr{Q}$) and the latter a notion of (perturbed) phase spectrum.

The next section focuses on the important case for applications in which the functions $\phi^{p,q}$ are built from correlation functions and for which these notions of power and phase spectra are further detailed through numerical computations and examples.

## VI. APPLICATION TO TIME SERIES: DATA-ADAPTIVE HARMONIC (DAH) SPECTRA

To simplify the presentation, we consider throughout this section a finite-dimensional Euclidean space $\mathscr{H}$.

### A. Correlation functions

In the case of a deterministic flow $(S_t)_{t\geq 0}$ acting on $\mathscr{H}$ and having an ergodic invariant measure $\mu$ supported by an attractor $\mathscr{A}$ [50, 51], let us recall that the correlation function $\rho^{f,g}(t)$ associated with two (sufficiently smoothed) observables $f, g : \mathscr{H} \to \mathbb{R}$, is given by

$$\rho^{f,g}(t) = \int_{\mathscr{H}} f(x)g(S_t x)\,\mathrm{d}\mu, \;\; t \geq 0, \, x \in \mathscr{H}. \tag{82}$$

In the case of an ergodic stochastic system, the correlation function is given by

$$\rho^{f,g}(t) = \int_{\mathscr{H}} f(x)P_t g(x)\,\mathrm{d}\mu, \;\; t \geq 0, \, x \in \mathscr{H}, \tag{83}$$

where $(P_t)_{t\geq 0}$ denotes the Markov semigroup possessing $\mu$ as (unique) ergodic invariant measure [52]. In each case, the function $\rho^{f,g}(t)$ is called a **correlation function** subordinated to the observables $f$ and $g$. Depending on these observables, various higher-order moment statistics can be considered.

In the deterministic setting, if $\mu$ is a Sinaï-Ruelle-Bowen measure [50, 51, 53], $\rho^{f,g}(t)$ corresponds to the more familiar cross-correlation coefficient at lag $t$ (for the observables $f$ and $g$) which takes then the equivalent form [54]

$$\rho^{f,g}(t) = \lim_{T\to\infty} \frac{1}{T} \int_0^T f(S_s x)g(S_{t+s}x)\,\mathrm{d}s, \tag{84}$$

for almost every $x$ in the basin attraction of $\mathscr{A}$ [51].

For an ergodic stochastic system, it takes the form

$$\rho^{f,g}(t) = \lim_{T\to\infty} \frac{1}{T} \int_0^T f(X_t^x)g(X_{t+s}^x)\,\mathrm{d}s, \tag{85}$$

$\mathbb{P}$-almost surely and for every $x$ in $\mathscr{H}$, when e.g. the Markov semigroup $P_t$ associated with the Markov process $X_t^x$ is strong Feller and irreducible; see [52].

### B. Periodized correlation functions, and associated semigroup kernels

We place ourselves in one of the working framework described in the previous subsection—i.e. either stochastic or deterministic—and apply the theory of Section IV to $\phi(t) = \rho^{f,g}(t)$ for a choice of two observables, $f$ and $g$. Note that if one chooses $f = g$, one recovers the familiar notion of lagged autocorrelation of e.g. the univariate time series $t \longmapsto f(S_t x)$. For this reason, we will sometimes, when $f \neq g$, name $\rho^{f,g}(t)$ as a cross-correlation function, to emphasize the difference with the notion of autocorrelation.

Obviously, the main assumption on which the theory of Sections III and IV relies on, is the periodicity of $\phi(t)$ which is not guaranteed in practice, except for periodic trajectories of course. However, motivated by physical applications, we are interested with mixing systems, i.e. systems for which a decay of correlations takes place. For stochastic (resp. deterministic) systems that are mixing, see e.g. [52] (resp. [51, 54]) and references therein.

For such systems, one should in practice choose $\tau$ large enough so that the decay of correlations has been sufficiently well-resolved. To apply the theory one needs then to periodize $\phi(t)$, and the resulting periodization $\varphi$ plays a role in the definition of the operator $\mathscr{L}_\varphi$. To understand it, let us consider an interval $I = [\alpha, \alpha+\tau]$ with $\alpha$ in $\mathbb{R}$. We define

$$\varphi(t) := \begin{cases} \phi(t), & t \in [\alpha, \alpha+\tau), \\ \phi(t-n\tau), & t \in [\alpha+n\tau, \alpha+(n+1)\tau),\ n \in \mathbb{Z}_*^+ \\ \phi(t+n\tau), & t \in (\alpha+(n-1)\tau, \alpha+n\tau],\ n \in \mathbb{Z}_*^-. \end{cases} \tag{86}$$

The left circular shift $T(t)$ defined in (31) applied to $\varphi$ gives thus the following expression in terms of $\phi$, for $r$ lying in the interval $I$

$$T(r)\varphi(s) := \begin{cases} \phi(r+s), & s \in [\alpha, \alpha+\tau-r), \\ \phi(r+s-\tau), & s \in [\alpha+\tau-r, \alpha+\tau). \end{cases} \tag{87}$$

In other words the kernel of the integral operator $\mathscr{L}_\varphi$ is modified compared with its expression given in (32).

As a choice of the base interval $I$ representing $\mathbb{T}$, a reasonable one consists of taking $\alpha = -\tau/2$ i.e.

$$I = \left[-\frac{\tau}{2}, \frac{\tau}{2}\right]. \tag{88}$$

This choice as a representation of $\mathbb{T}$ relies in part on the symmetry of $\phi(t) = \rho^{f,g}(t)$, i.e. $\rho^{f,g}(t) = \rho^{f,g}(-t)$. Thus $\phi(-\tau/2) = \phi(\tau/2)$. Other choices, like $\alpha = 0$, results—due to the assumed decay of correlations—with "jumps" that occur at the junctions points in the periodization of $\phi$ given in (86).

These jumps do not constitute however an obstruction to apply the theory. Theorem IV.1 allows indeed for a discontinuity of this type for $\varphi$ used in the definition of $\mathscr{L}_\varphi$. The reason is tied to the characterization of $D(A)$ given in (30) (with$g \equiv 1$) that shows that a function $\varphi$ with a discontinuity is allowed, as long as $\varphi$ defined in (86) is absolutely continuous and $\varphi'$ lies in $L^2(\mathbb{T})$.

In Sect. VI D below, we consider the choice of the base interval $I$ given by (88) in the case of a multidimensional DAH spectrum obtained from a collection of correlation functions. We turn next to the case $I = [0, \tau]$, i.e. $\alpha = 0$, for which the theory allows to recover instructive results from the literature on Hankel matrices. In that respect we address first the discretization of the operator $\mathscr{L}_\varphi$.

### C. DAH spectrum from a single correlation function: discretization

Given a uniform partition $t_0 = 0 < \cdots < t_{N-1} < \tau = t_N$, we consider the family of functions

$$\chi_j(s) := \begin{cases} 1, & s \in [t_j, t_{j+1}), \\ 0, & \text{else}. \end{cases} \quad 0 \le j \le N-1. \tag{89}$$

We simply denote by $\rho$, the correlation function $\rho^{f,g}$. By the definition (87) of the left circular shift acting on the periodization $\varphi$ in (86) (with $\alpha = 0$), we have from (32),

$$\mathscr{L}_\varphi(\Psi)(r) = \frac{1}{\tau}\left(\int_0^{\tau-r} \rho(s+r)\Psi(s)\,\mathrm{d}s + \int_{\tau-r}^{\tau} \rho(r+s-\tau)\Psi(s)\,\mathrm{d}s\right). \tag{90}$$

Now let us approximate $\rho$ by a piecewise constant function $f$, as follows

$$f(s) = \rho(t_j), \text{ if } s \in [t_j, t_{j+1}). \tag{91}$$

Take $t_j = j\delta\tau$ with $\delta\tau = \tau/N$, then

$$\mathscr{L}_\varphi[\chi_j](r) \approx \frac{\delta\tau}{\tau}\rho((j+k)\delta\tau) \mod (\tau), \text{ if } r \in [t_k, t_{k+1}). \tag{92}$$

Denoting $\rho(l\delta)$ by $\rho_l$, the operator $\mathscr{L}_\varphi$ can be thus approximated by the following Hankel $N \times N$ matrix

$$\mathbf{H} = \frac{1}{N}\begin{pmatrix} \rho_0 & \rho_1 & \rho_2 & \cdots & \rho_{N-1} \\ \rho_1 & \rho_2 & \rho_3 & \iddots & \rho_0 \\ \rho_2 & \vdots & \iddots & \iddots & \rho_1 \\ \vdots & \iddots & \iddots & \iddots & \vdots \\ \rho_{N-1} & \rho_0 & \rho_1 & \dots & \rho_{N-2} \end{pmatrix}. \tag{93}$$

Without any surprise, $\mathbf{H}$ is a left circulant matrix, that we denote sometimes as

$$\mathbf{H} = l\text{-circ}(\rho_0, \rho_1, \cdots, \rho_{N-1}), \tag{94}$$

or in other words, the rows of $\mathbf{H}$ are obtained by successive shifts to the left of one position starting from the row $r = (\rho_0, \rho_1, \cdots, \rho_{N-1})$ as a first row.

Let $\mathbf{T}$ be the right circulant matrix which has the same first row than the left circulant matrix $\mathbf{H}$. Then $\mathbf{T}$ is Toeplitz and $\mathbf{H} = P\mathbf{T}$ for the permutation matrix $P = 1 \oplus J_{n-1}$, where $J_{N-1}$ is the reversal $(N-1)\times(N-1)$ matrix obtained by flipping the $(N-1)\times(N-1)$ identity matrix $I_{N-1}$ from left to right. Hence $\mathbf{H}$ and $\mathbf{T}$ have identical singular values. However, circulant matrices and real anticirculant matrices (which are also real symmetric) are normal matrices, so they can be unitarily diagonalized and their singular values are the moduli of their eigenvalues. Therefore, the eigenvalues $\mathbf{H}$ and $\mathbf{T}$ have identical moduli. Finally, since real anticirculant matrices are real symmetric, $\mathbf{H}$ has real eigenvalues.

The spectral theory of circulant matrices is well established although somehow scattered in the literature; see e.g. [55–57]. The inherent periodicity of circulant matrices relates them to Fourier analysis and group theory. In contrast to the aforementioned standard circulant Toeplitz matrix $\mathbf{T}$ arising in Fourier analysis (see e.g. [58, Def. 5.12]), the left-circulant matrix $\mathbf{H}$ is not diagonalized by the unitary Fourier matrix [57, Theorem 3.5] and, as Proposition VI.1 summarizes, its eigenvectors depend actually on the data and more exactly relate to the phase spectrum as a counterpart of Theorem IV.1 and a consequence of formula (43).

Let us introduce the Fourier frequency $\omega_k = 2\pi k/N$, and the periodogram of $\{\rho_j\}$, namely

$$I_N(\omega_k) := \frac{1}{N}\Big|\sum_{j=0}^{N-1} \rho_j e^{-ij\omega_k}\Big|^2, \;\; 0 \leq k \leq N-1, \tag{95}$$

with $i^2 = -1$.

We have then

**Proposition VI.1** *The eigenvalues $\{\lambda_k\}_{k\in\{0,\ldots,N-1\}}$ of the Hankel matrix $\mathbf{H}$ given by (93) can be arranged in terms of Fourier frequencies so that*

$$\begin{aligned} \lambda_0 &= N^{-1}\sum_{n=0}^{N-1}\rho_n, \\ \lambda_{\frac{N}{2}} &= N^{-1}\sum_{n=0}^{N-1}(-1)^n\rho_n, \; \textit{if } N \textit{ is even}, \end{aligned} \tag{96}$$

*and for $1 \leq k \leq \left[\frac{N-1}{2}\right]$,*

$$\lambda_k = -\lambda_{N-k} = N^{-1/2}\sqrt{I_N(\omega_k)}, \tag{97}$$

*where $[x]$ denotes the largest integer less than or equal to $x$.*

*Furthermore, for each pair $(\lambda_k, \lambda_{N-k})$ corresponding thus to a Fourier frequency $\omega_k = 2\pi k/N$, the corresponding pair of eigenvectors $(v_k, v_{N-k})$ satisfies*

$$\begin{aligned} &\arg(\widehat{\rho}(\omega_k)) = \arg(\lambda_j\widehat{v_j}(\omega_k)) - \arg(\widehat{v_j}(\omega_k)), \\ &\quad \textit{with } j \in \{k, N-k\}, \; 1 \leq k \leq \left[\frac{N-1}{2}\right], \end{aligned} \tag{98}$$

*where $\widehat{u}(\omega)$ denotes the discrete Fourier transform—assessed at the frequency $\omega$—of a vector $u$.*

**Proof.** For a proof of (96)-(97) we refer to [56, Lemma 1]. The identity (98) is the analogue of (51) in the discrete setting. ∎

In what follows the eigenvalues of the matrix $\mathbf{H}$ when ranked from the lowest to the highest resolved frequency, form what we call the **DAH power spectrum**. The **DAH phase spectrum** is formed by the phase of the corresponding eigenvectors, also ranked per frequency. Whereas, due to Proposition VI.1, these concepts coincide for one-dimensional signals with the standard notions of power and phase spectra (see also Remark VI.2-(ii) below), their multidimensional generalization following Sect. V and applied to a collection of cross-correlation functions, lead to more subtle but still useful interpretations. This generalization is discussed in in further details in Section VI D below.

**Remark VI.1** *Let us introduce the standard cyclic permutation matrix,*

$$\mathbf{K} := \begin{pmatrix} 0 & 1 & \cdots & 0 & 0 \\ 0 & 0 & \ddots & \vdots & 0 \\ 0 & 0 & \ddots & 1 & \vdots \\ \vdots & \ddots & \ddots & 0 & 1 \\ 1 & \cdots & 0 & 0 & 0 \end{pmatrix}. \tag{99}$$

*Then $\mathbf{H}$ is naturally associated with the Krylov space $\mathscr{K}(\mathbf{K},\rho)$ generated by the matrix $\mathbf{K}$ given by (99) and the vector $\rho = (\rho_0, \cdots, \rho_{N-1})^{\mathrm{tr}}$, namely*

$$\mathscr{K}(\mathbf{K},\rho) = \mathrm{span}(\rho, \mathbf{K}\rho, \mathbf{K}^2\rho, \cdots, \mathbf{K}^{\mathrm{N}-1}\rho). \tag{100}$$

*Indeed, the $j$-th column of $\mathbf{H}$ is exactly given by $\mathbf{K}^{j-1}\rho$.*

**Remark VI.2**

*(i) Proposition VI.1 remains valid when other base intervals $I$ than $[0,\tau]$ are used. The modification consists essentially of replacing accordingly the summation bounds in the definition of the periodogram $I_N$ given by (95).*

*(ii) Let $t \mapsto X_t^x$ denote a stochastic process emanating from $x$ in $\mathscr{H} := \mathbb{R}^q$ $(q > 1)$ and generated by a system of SDEs. Given two observables $f, g : \mathscr{H} \to \mathbb{R}$, and recalling that $\rho$ denotes the correlation function $\rho^{f,g}$ given by (83), we infer from the formula (97) that the $\lambda_k$'s provide the* **(one-sided) cross-power spectrum**, *$\Gamma_{f,g}$, of the real-valued stochastic processes $t \mapsto f(X_t^x)$ and $t \mapsto g(X_t^x)$.*

*The DAH eigenvectors provide furthermore the* **(one-sided) cross phase spectrum** *via the formula (98), in the sense that*

$$\begin{aligned} \Gamma_{f,g}(\omega_k) = \sqrt{N}|\lambda_j| \exp\Big(\arg(\lambda_j\widehat{v_j}(\omega_k)) - \arg(\widehat{v_j}(\omega_k))\Big), \\ \textit{with } j \in \{k, N-k\}, \; 1 \leq k \leq \left[\frac{N-1}{2}\right]. \end{aligned} \tag{101}$$

### D. Multidimensional DAH spectrum from a collection of correlation functions

Given a base interval $I$ to be $[\alpha, \alpha+\tau]$, one considers the periodization $\varphi^{f,g} = \varphi$ given in (86) for $\phi = \rho^{f,g}$, with $\rho^{f,g}$ denoting a correlation function given either by (84) or by (85), depending on the context.

The operator $\mathscr{L}_\phi$ defined in (32) becomes then

$$\begin{aligned} \mathscr{L}_{\varphi^{f,g}}(\Psi)(r) = &\frac{1}{\tau}\Bigg(\int_\alpha^{\alpha+\tau-r} \rho^{f,g}(s+r)\Psi(s)\,\mathrm{d}s \\ &+ \int_{\alpha+\tau-r}^{\alpha+\tau} \rho^{f,g}(r+s-\tau)\Psi(s)\,\mathrm{d}s\Bigg), \end{aligned} \tag{102}$$

where $\Psi$ is any arbitrary function of $L^2(\mathbb{T})$.

By choosing $I = [-\tau/2, \tau/2]$ (i.e. $\alpha = -\tau/2$), a uniform partition

$$-\frac{\tau}{2} = t_{-M+1} \leq \cdots \leq t_{-1} \leq t_0 = 0, \qquad 0 \leq t_1 \leq \cdots \leq t_{M-1} \leq \frac{\tau}{2} = t_M, \tag{103}$$

and elementary functions such as in (89), each operator $\mathscr{L}_{\varphi^{i,j}}$ is approximated by the following Hankel $(2M-1)\times(2M-1)$ matrix

$$\mathbf{H}^{(i,j)} = \begin{pmatrix} \rho^{(i,j)}_{-M+1} & \rho^{(i,j)}_{-M+2} & \cdots & \rho^{(i,j)}_{0} & \rho^{(i,j)}_{1} & \cdots & \rho^{(i,j)}_{M-1} \\ \rho^{(i,j)}_{-M+2} & \iddots & \iddots & \iddots & \iddots & \iddots & \rho^{(i,j)}_{-M+1} \\ \vdots & \iddots & \iddots & \iddots & \iddots & \iddots & \rho^{(i,j)}_{-M+2} \\ \rho^{(i,j)}_{0} & \iddots & \iddots & \iddots & \rho^{(i,j)}_{-M+1} & \iddots & \vdots \\ \rho^{(i,j)}_{1} & \iddots & \iddots & \iddots & \rho^{(i,j)}_{-M+2} & \iddots & \rho^{(i,j)}_{0} \\ \vdots & \rho^{(i,j)}_{M-1} & \rho^{(i,j)}_{-M+1} & \iddots & \iddots & \iddots & \vdots \\ \rho^{(i,j)}_{M-1} & \rho^{(i,j)}_{-M+1} & \rho^{(i,j)}_{-M+2} & \cdots & \rho^{(i,j)}_{0} & \cdots & \rho^{(i,j)}_{M-2} \end{pmatrix}. \tag{104}$$

The operator $\mathfrak{L}$ defined in (53)-(54)—in which each operator $\mathscr{L}_{\phi^{p,q}}$ therein is replaced by $\mathscr{L}_{\varphi^{p,q}}$ given in (102) above—is then approximated by the following block-Hankel matrix $\mathfrak{C}$ formed by $d^2$ blocks, $\mathfrak{C}^{(i,j)}$, each of size $(2M-1)\times(2M-1)$ and given as follows

$$\mathfrak{C}^{(i,j)} = \mathbf{H}^{(i,j)}, \text{ if } i \leq j, \qquad \mathfrak{C}^{(i,j)} = \mathbf{H}^{(j,i)}, \text{ if } j < i, \tag{105}$$

for $(i,j)$ in $\{1,\cdots,d\}^2$.

Hereafter we use $M' = 2M-1$ for concision, reindexing possibly from 1 to $M'$ the string $\{-M+1,\cdots,M-1\}$. By construction and in agreement with Lemma V.1, the matrix $\mathfrak{C}$ is real symmetric, not circulant with each block that is Hankel.

Since the matrix $\mathfrak{C}$ is the discretization of the operator $\mathfrak{L}$ (defined in (53)-(54) and in which each operator $\mathscr{L}_{\phi^{p,q}}$ therein is replaced by $\mathscr{L}_{\varphi^{p,q}}$ given in (102)), Theorem V.1 is enlightening concerning the spectrum of $\mathfrak{C}$. In particular the eigenvalues of $\mathfrak{C}$ relate according to (59) to the singular values of the **symmetrized cross-spectral matrix** $\Phi_\delta(f_k)$ whose entries are still given by (58) except that $\widehat{\phi^{p,q}}(f_k)$ is obtained here as the **discrete** Fourier transform (evaluated at the frequency $f_k$) of the correlation function $\rho^{p,q}(t)$, after its periodisation following Sect. VI B.

As noted earlier for the operator $\mathfrak{L}$ and applied here to the matrix $\mathfrak{C}$, Theorem V.1 and Remark V.1-(ii) ensure that the **eigenvalues** of $\mathfrak{C}$ relate naturally to the **quadratic form** $\mathscr{Q}_\delta$ given by

$$\mathscr{Q}_\delta(x) := \mathrm{Re}(x^T \Phi_\delta(f_k) x), \;\; x \in \mathbb{C}^d, \tag{106}$$

via the variational principle given in (76)-(77), in which $\Phi_\delta(f_k)$ replaces $\Phi(f_k)$. The eigenvalues of $\mathfrak{C}$ convey thus also an idea of energy distribution per frequency, as measured with the quadratic form $\mathscr{Q}_\delta$.

Grouping per frequency the absolute value of the eigenvalues of the matrix $\mathfrak{C}$, and displaying those eigenvalues from the lowest to the highest resolved frequency, we obtain what we call the **multidimensional DAH power spectrum**, omitting the adjective "multidimensional" when this is clear from the context.

The **multidimensional DAH phase spectrum** is obtained by grouping per frequency the $\eta^p_+(f_k)$ (resp. $\eta^p_-(f_k)$) as discussed at the end of Sect. V. Note that in practice, the $\eta^p_+(f_k)$ (resp. $\eta^p_-(f_k)$) are calculated by evaluating the LHS of (80) in which $E^k_p$ is replaced by the $p$-th "snippet" of size $d$, extracted from an eigenvector (still denoted by) $W_k$ of $\mathfrak{C}$. In applications, the eigenvectors of $\mathfrak{C}$ are called the **DAH modes (DAHMs)**, and the $E^k_p$ are called **DAHM snippets**; see also Appendix A.

Given a $d$-dimensional time series $\mathbf{X}(t)$ from which the matrix $\mathfrak{C}$ is formed, we will also pay attention to another important quantity, namely the **DAH coefficients (DAHCs)** obtained by projection of $\mathbf{X}(t)$ onto the DAHMs; see Appendix A for more details. These DAHCs play an important role in Section VII concerned with inverse modeling.

For the moment, we turn next to basic examples in order to illustrates the notions of power and phase spectra, as well DAHMs and DAHCs to help build up first intuitions.

### E. DAH spectrum of traveling waves with modulated wavelength

We consider the following one-dimensional time-periodic scalar field, given for $x$ in $[0,L]$ and $t > 0$ by

$$q(x,t) = \mathrm{Re}\Big(\exp\big(i\big(k(x)x - \omega_c t\big)\big)\Big), \tag{107}$$

with a wave length $k(x)$ depending on the spatial location. Given a uniform discretization $[0,L]$ by $N_x$ points for a mesh size $\delta x$, we choose hereafter $k(x)$ to be given according to

$$k(x) = \frac{x}{\delta x} - p\alpha, \;\; \alpha = \frac{N_x}{m}, \;\; 1 \leq m \leq N_x. \tag{108}$$

Here $p$ denotes the integer for which the quotient $x/(\alpha\delta x)$ is within the roundoff error of that integer. The parameter $\omega_c$ is the Fourier frequency given by

$$\omega_c = 2\pi\frac{(l-1)}{K'},$$

with $K' = 2K-1$ and $1 \leq l \leq M'$. For the experiments reported below, we set $L = 5\pi$, $N_x = 42$, $K = 15$, $l = 3$ and $m = 19$. For simplicity, a time-stepping, $\delta t = 1$, is chosen to discretize the time-variable $t$ in (107). The resulting patterns are shown in Fig. 1-(a). As a benchmark, a traveling wave with $k(x)$ replaced by its mean value $\bar{k}$, is shown in Fig. 1-(b). The effect of $k(x)$ is thus to distort a spatial cosine wave, $\cos(\bar{k}x - \omega_c t)$, by introducing some local oscillations occurring at scales smaller than $1/\bar{k}$.

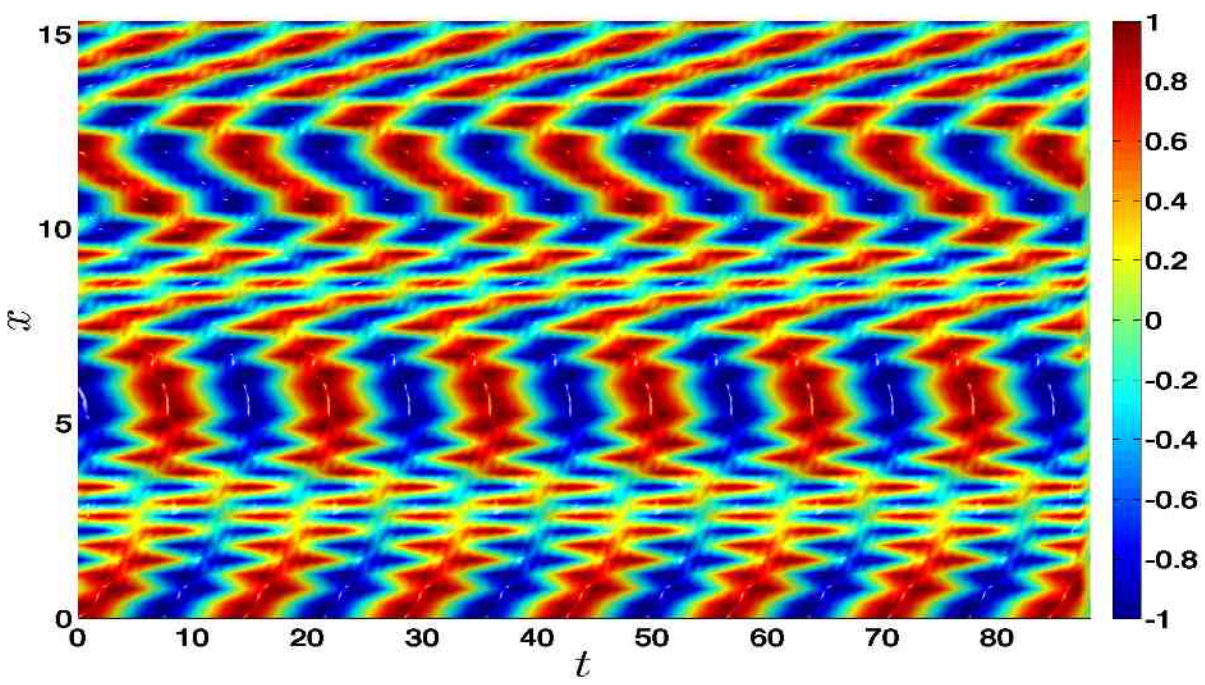

(a) **Modulated wave**

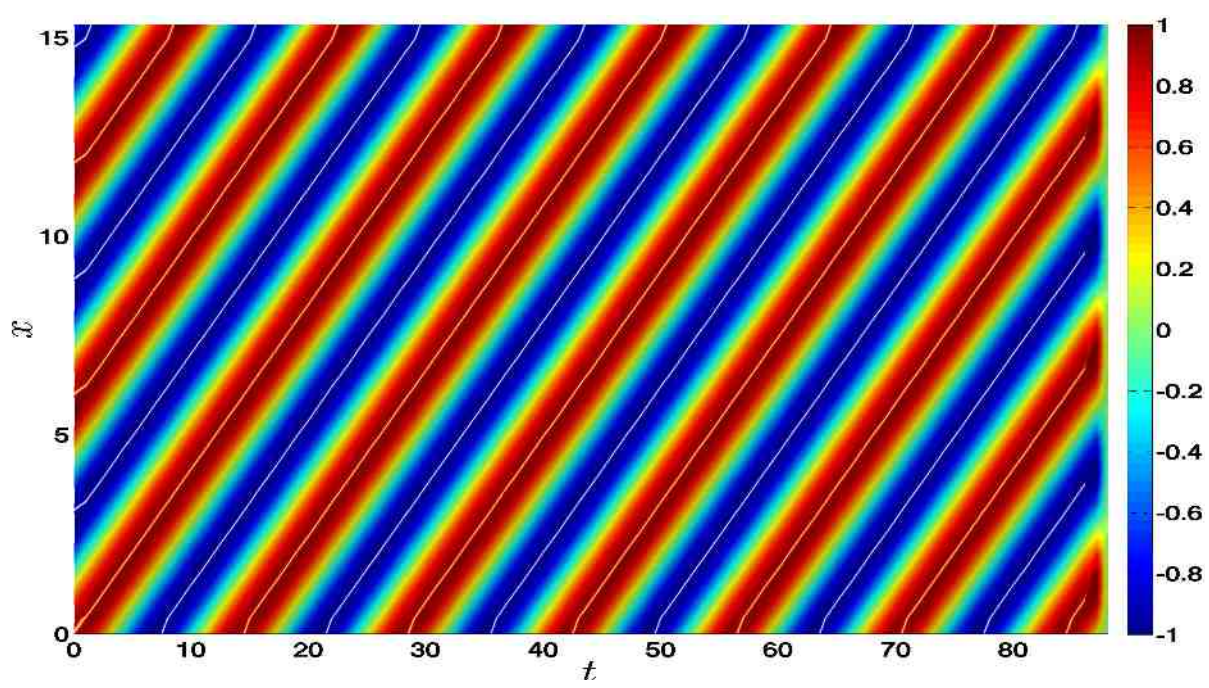

(b) **Travelling wave**

FIG. 1: **Spatio-temporal fields.** Panel (a) shows a time evolution of $q(x,t)$ with modulated wave lengths $k(x)$ given by (108). Panel (b) shows that obtained replacing $k(x)$ by its mean value $\bar{k}$; see text for more details.

Compared to the case of a pure traveling wave shown in Fig. 1-(b), some channels in Fig. 1-(a) still exhibit collectively some traveling wave patterns (between e.g. $x \approx 14$ and $x = 5\pi$), while others exhibit now standing-like waves patterns as it can be observed between $x \approx 4$ and $x \approx 7$.

The corresponding DAH power spectrum and DAH phase spectrum are shown in Fig. 2-(a) and Fig. 2-(b), respectively. Their computations has been obtained by diagonalizing the matrix $\mathfrak{C}$ given in (105) and following the procedure described in Sect. VI D, in which the cross correlations are estimated from $q(x_j,t_n)$ obtained from (107), with $x_j = (j-1)L/N_x$, $1 \leq j \leq N_x$, and $t_n = n$, for $0 \leq n \leq 3\times 10^3$. The corresponding DAHMs and DAHCs are shown in Fig. 4.

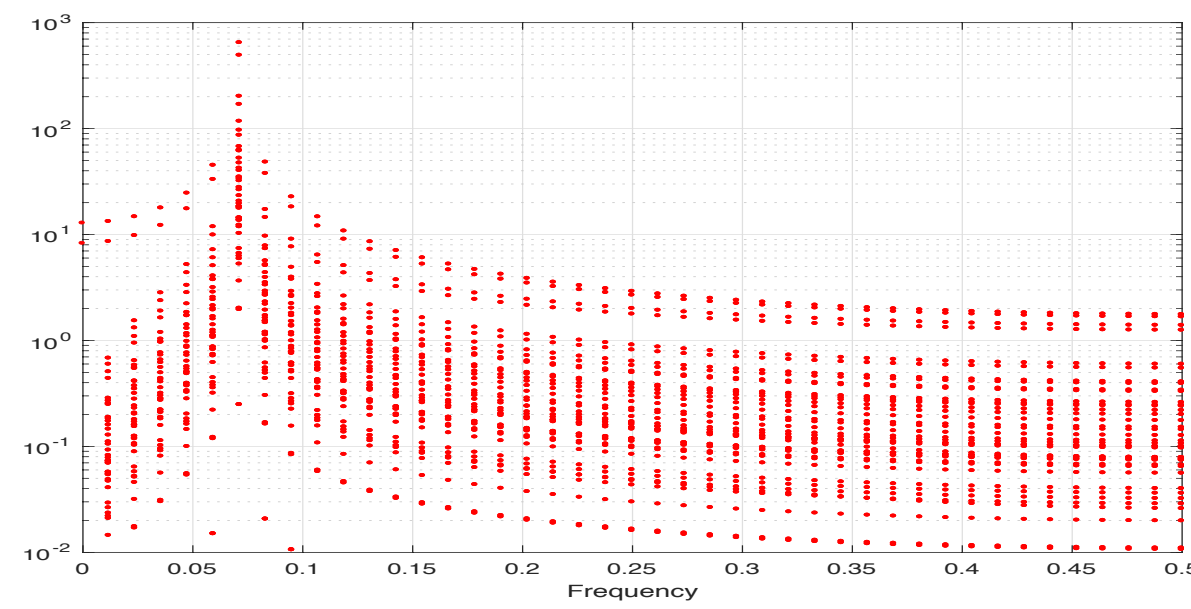

(a) **DAH power spectrum**

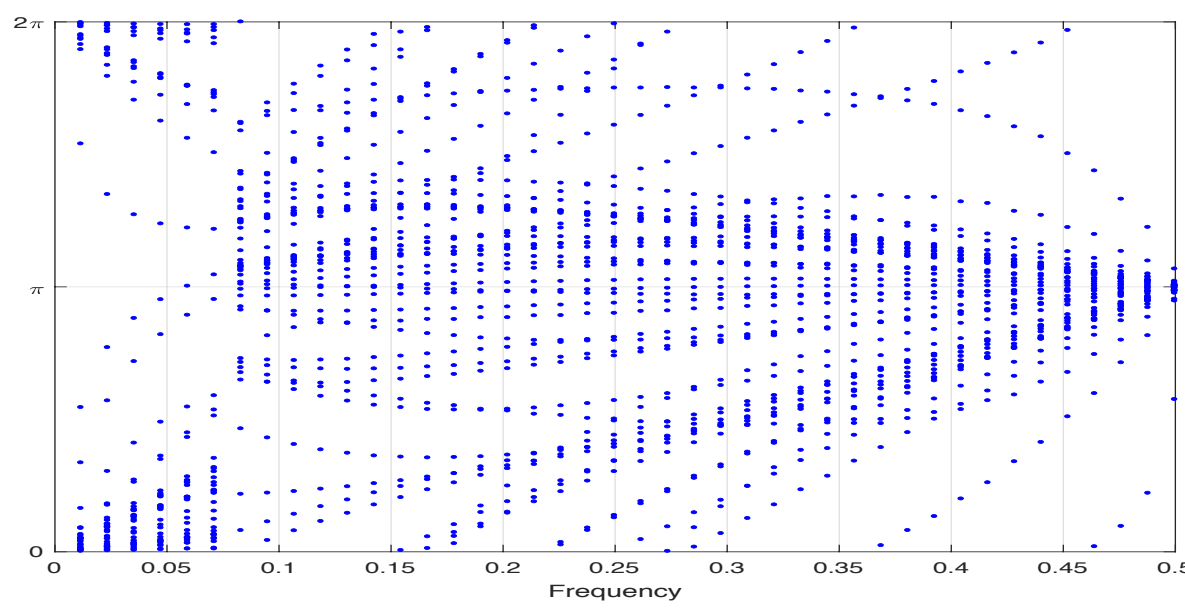

(b) **DAH phase spectrum**

FIG. 2: **Multidimensional DAH power and phase spectra.** The DAH spectra shown here are those associated with the modulated wave shown in Fig. 1-(a). We chose $M = 43$ to estimate the cross correlations and form the corresponding matrix $\mathfrak{C}$ given in (105). Each non-zero frequency is associated with 42 eigenpairs that correspond to the total of $d = N_x = 42$ channels used to generate the field $q(x,t)$ from (107).

### 1. *DAH power spectrum and DAHMs*

We first comment on the DAH power spectrum and DAHMs. The main noticeable feature of the DAH power spectrum shown in Fig. 2-(a) (resp. Fig. 3-(a)) is the presence of a peak located at the (temporal) frequency close to $f_c = \omega_c/(2\pi)$ (due to our discretization), the characteristic frequency of the spatio-temporal field $q(x,t)$ given by (107) (resp. (107) with $k(x)$ replaced by $\bar{k}$). Recalling the variational characterization of the eigenvalues of $\mathfrak{C}$, one has thus here an illustration of the relevance of the energy $\mathscr{Q}_\delta$ given by (106) for indicating the presence of a dominant oscillation; the latter being expressed by the presence of a dominant pair of eigenvalues (i.e. with the highest max-min for $\mathscr{Q}_\delta$; see Remark V.1-(ii)) at that frequency.

This dominant pair of eigenvalues come with an energetic oscillatory pair of DAHCs as shown in the panel located at the intersection of the bottom row and third column of Fig. 4 (resp. Fig. 5) for $f = 0.071 \approx \omega_c/(2\pi)$; compare with the DAHCs of the third column for other frequencies.

For the modulated-wave case, the mixture of traveling/standing-like wave patterns as well as the local small-scale oscillations are outstandingly well captured, at each frequency $f \neq 0$, by the pair of DAHMs associated with the pair of dominant eigenvalues; see first and second column of Fig. 4. We emphasize indeed that this mixture of patterns present in the field $q(x,t)$ shown in Fig. 1-(a), is not only found for the (few) mode(s) associated with the DAHCs of largest variance(s), but rather is distributed across a broad frequency band of Fourier frequencies $f$.

Noticeably, these patterns scale with the frequency: at each frequency $f$, the pair of dominant DAHMs exhibit the same patterns that repeat with a $1/f$-temporal periodicity; see 1st and 2nd columns of Fig. 4. The same elemental spatial pattern is thus repeated for a pair of (dominant) DAHMs across the frequencies.

The same feature is observed in the case of a pure traveling wave shown in Fig. 1-(b). Here again, the (dominant) DAHMs capture the dominant pattern displayed by the field $q(x,t)$ across a broad range of frequencies; see 1st and 2nd columns of Fig. 5.

Such a scaling property (in time)—of the principal pattern contained in the original dataset—as revealed by the dominant pairs of DAHM across the frequencies, is in sharp contrast with other modes that would be obtained by techniques relying on a decomposition of the variance.

Such methods applied to $q(x,t)$ would obviously identify another type of leading modes that would capture most of the variance, without however displaying a scaled-version of the dominant pattern across the frequencies, in particular for the modes that capture a small fraction of the variance. The third columns of Fig. 4 and Fig. 5 show to the contrary that even for pairs of DAHMs associated with DAHCs of small amplitude, these modes exhibit the dominant pattern, repeated at the associated frequency $f$; see the cases $f = 0.024$ and $f = 0.048$ in the 1st and 2nd columns of Fig. 4 or Fig. 5.

#### 2. *DAH phase spectrum*

We comment now the DAH phase spectra shown in Fig. 2-(b) and Fig. 3-(b). In both cases, the most striking feature about these phase spectra is the distinction appearing in their shape by comparing the frequency band $0 < f < f_c$ with $f > f_c$. For $f > f_c$ one can indeed distinguish that all the lines converge towards a single point as one approaches the Nyquist frequency $f_N$.

By relying on the discussion conducted at the end of Section V, this convergence corresponds to a convergence of the $\eta_+^p(f)$ towards a same value, expressing thus a form of **phase synchronization** between the channels of a DAHM as $f \to f_N$, by recalling (81) and (61).

One could think that such a phase synchronization revealed by the DAH phase spectrum is limited to the case of pure or modulated traveling waves, but actually, several conducted experiments on other spatio-temporal fields issued from various dynamical systems (not shown), seem to indicate that this feature occurs for a broad set of deterministic dynamics, including chaotic ones; see also Section VIII A 2 below. Roughly speaking, the poor resolution that comes with the Nyquist frequency explains that differences between the phases across the channels—subtle by nature—must vanish as $f \to f_N$. However, here lies maybe the most useful aspect of the DAH phase spectrum. As we will see in Sect. VIII B 2, dataset issued from a stochastic system does not share this convergence property of the phases across the channels of a DAHM, as $f \to f_N$. The DAH phase spectrum appears thus to be a useful diagnostic tool to distinguish between spatio-temporal deterministic or stochastic dynamics.

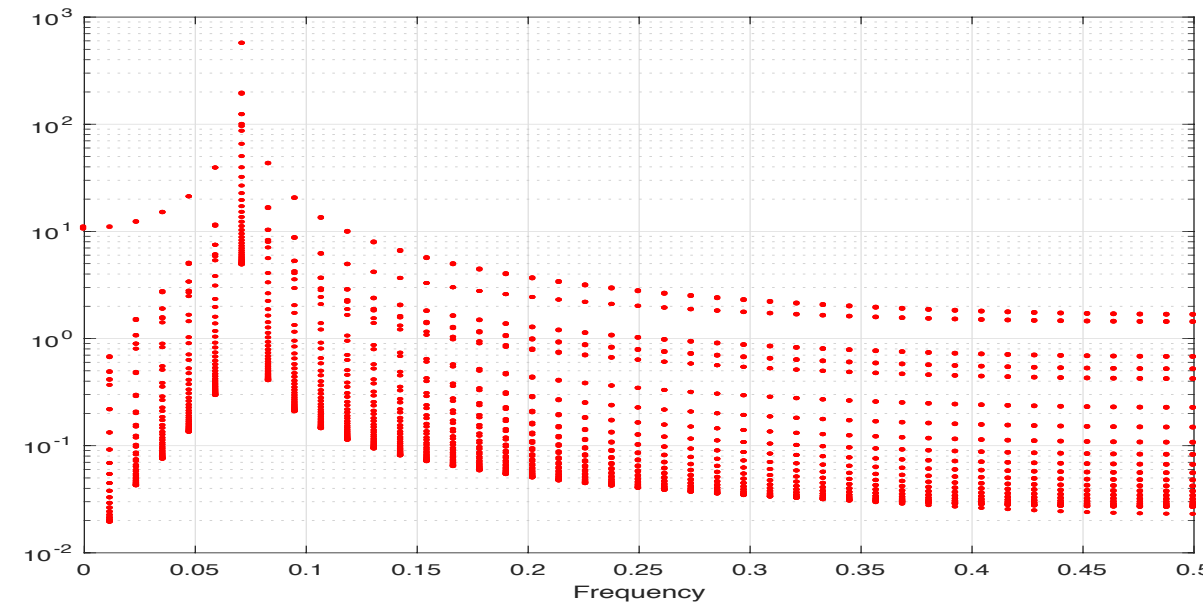

(a) **DAH power spectrum**

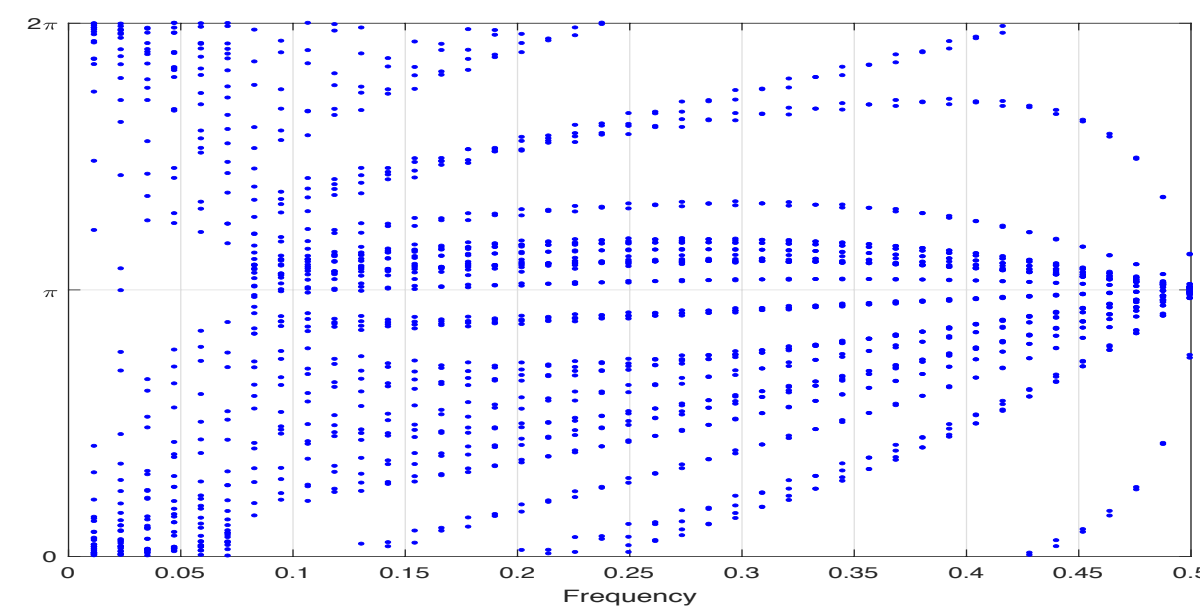

(b) **DAH phase spectrum**

FIG. 3: **Multidimensional DAH power and phase spectra.** The DAH spectra shown here are those associated with the pure traveling wave shown in Fig. 1-(b). We chose $M = 43$ to estimate the cross correlations and form the corresponding matrix $\mathfrak{C}$ given in (105). The DAH power spectrum and DAH phase spectrum are obtained as described in Sect. VI D.

## VII. INVERSE STOCHASTIC MULTILAYER STUART-LANDAU MODELS

### A. Preliminaries

Let $\mathbf{X}(t) = (X_1(t), \cdots, X_d(t))$ be a $d$-dimensional time series from which the operator $\mathfrak{L}$ defined in (53)-(54) is built from the cross-correlation functions $\rho^{i,j}(t)$, after (possible) periodization such as described in Sect. VI B. We adopt the notations of Theorem V.1 and form the following time-continuous version of the DAHC associated with an eigenfunction $W_k$ (see (A1) in Appendix A), namely

$$\begin{aligned}\xi^k(t) &:= \sum_{p=1}^{d} \int_{\mathbb{T}} X_p(t+s) E_p^k(s)\,\mathrm{d}s \\ &= \sum_{p=1}^{d} \int_{\mathbb{T}} X_p(t+s) B_p^k \cos(\omega_k s + \theta_p^k)\,\mathrm{d}s \\ &= \sum_{p=1}^{d} \int_{\mathbb{T}} T(t) X_p(s) B_p^k \cos(\omega_k s + \theta_p^k)\,\mathrm{d}s.\end{aligned} \tag{109}$$

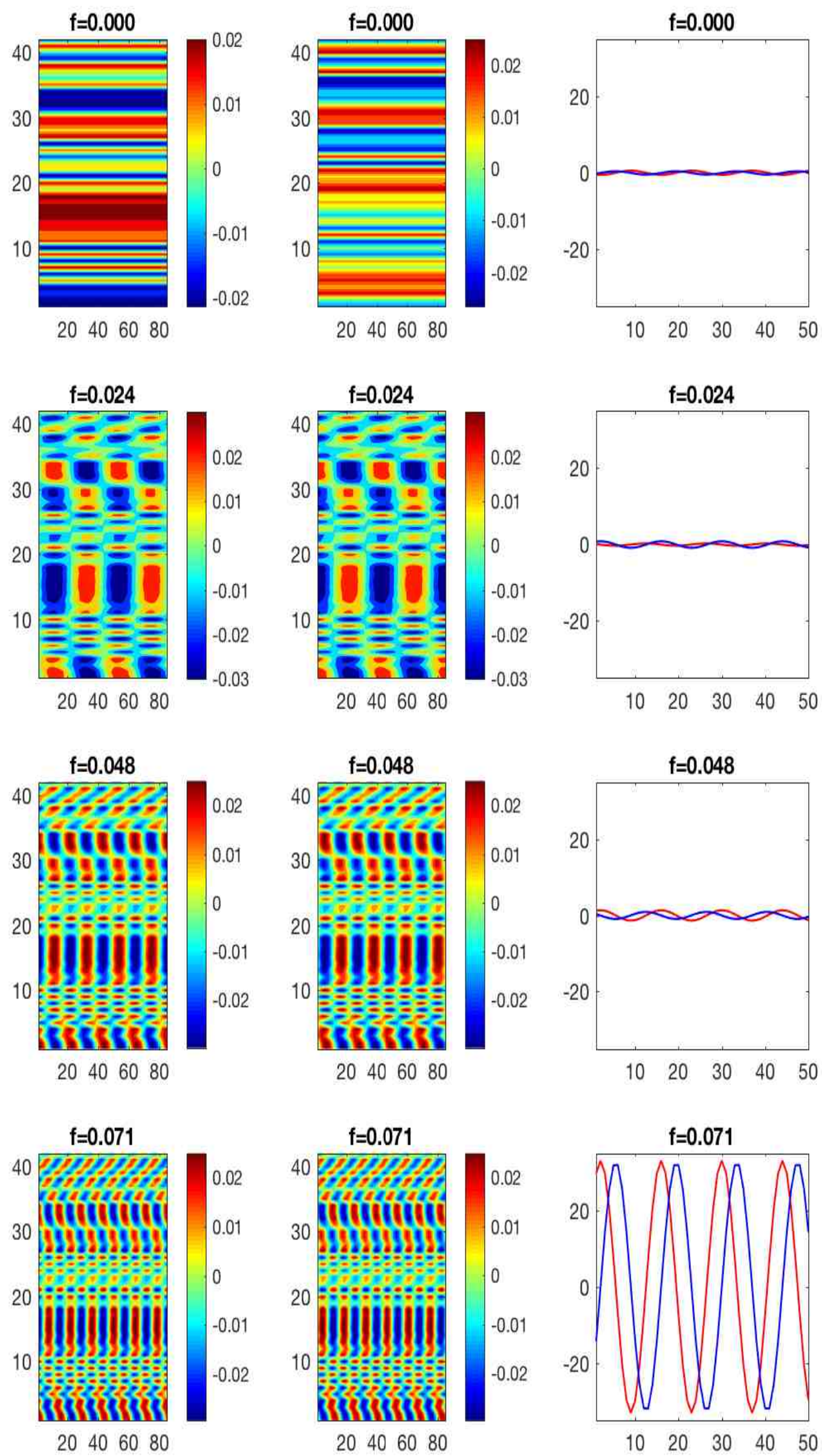

FIG. 4: **DAHMs and DAHCs from the modulated wave of Fig. 1-(a).** Each row of the 1st and 2nd column, shows a pair of DAHMs at a given frequency such as indicated. There, the horizontal axes represent the temporal embedding dimension $M' = 2M - 1$ (here $M = 43$), while the vertical ones indicate the spatial channel $1 \leq k \leq d$; $d = N_x = 42$. The DAHMs and DAHCs shown here correspond, for each indicated frequency $f$, to the dominant pair of DAH eigenvalues located above $f$ in Fig. 2-(a).

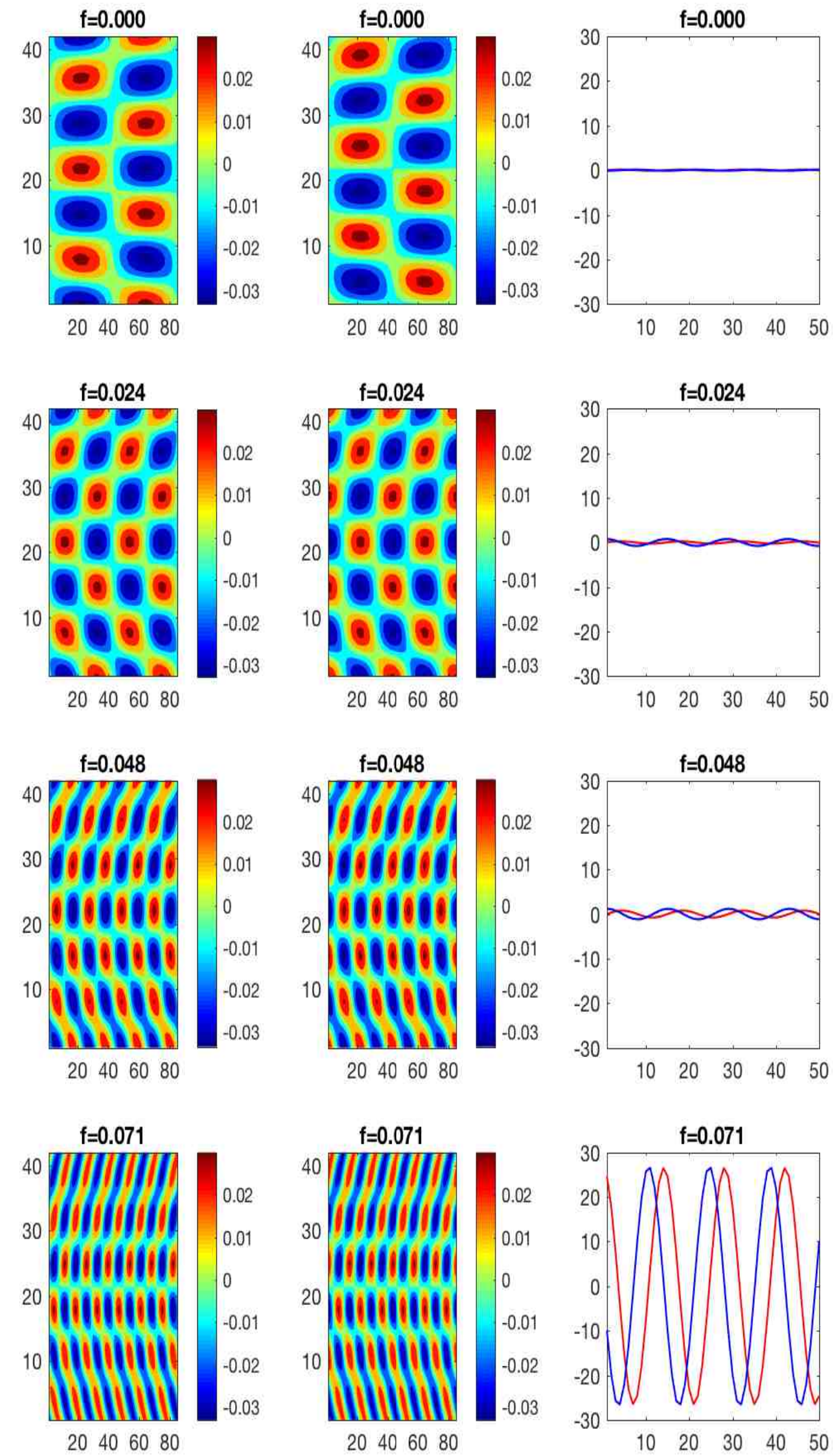

FIG. 5: **DAHMs and DAHCs from the traveling wave of Fig. 1-(b).** Same as in Fig. 4 but for the traveling wave of Fig. 1-(b). The DAHMs and DAHCs shown here correspond, for each indicated frequency $f$, to the dominant pair of DAH eigenvalues located above $f$ in Fig. 3-(a).

If $X_p$ lies in $D(A)$ given in (30) (with $g = 1$), then

$$\frac{\mathrm{d}T(t)X_p}{\mathrm{d}t} = AT(t)X_p = T(t)AX_p, \tag{110}$$

see e.g. [35, Chap. II, Lemma 1.3].

In particular

$$\frac{\mathrm{d}\xi^k(t)}{\mathrm{d}t} = \sum_{p=1}^{d} \int_{\mathbb{T}} T(t) \frac{\mathrm{d}X_p}{\mathrm{d}s} B_p^k \cos(\omega_k s + \theta_p^k)\,\mathrm{d}s. \tag{111}$$

Similarly,

$$\frac{\mathrm{d}\widetilde{\xi^k}(t)}{\mathrm{d}t} = \sum_{p=1}^{d} \int_{\mathbb{T}} T(t) \frac{\mathrm{d}X_p}{\mathrm{d}s} B_p^k \sin(\omega_k s + \theta_p^k)\,\mathrm{d}s \tag{112}$$

Thus

$$\frac{\mathrm{d}z^k(t)}{\mathrm{d}t} = \sum_{p=1}^{d} \int_{\mathbb{T}} T(t) \frac{\mathrm{d}X_p}{\mathrm{d}s} B_p^k \exp(i(\omega_k s + \theta_p^k))\,\mathrm{d}s \tag{113}$$

An integration by parts gives

$$\frac{\mathrm{d}z^k(t)}{\mathrm{d}t} = -i\omega_k \sum_{p=1}^{d} \int_{\mathbb{T}} T(t) X_p(s) B_p^k \exp(i(\omega_k s + \theta_p^k))\,\mathrm{d}s, \tag{114}$$

that is

$$\frac{\mathrm{d}z^k(t)}{\mathrm{d}t} = -i\omega_k z^k, \tag{115}$$

supplemented by initial condition

$$z^k(0) = \sum_{p=1}^{d} \int_{\mathbb{T}} X_p(s) B_p^k \exp(i(\omega_k s + \theta_p^k))\,\mathrm{d}s. \tag{116}$$

If $X_p$ is not in the domain $D(A)$ (such as for $\mathbf{X}(t)$ issued from a mixing dynamical system) and the base interval $I$ for $\mathbb{T}$ is taken to be $[-\tau, 0]$, then boundaries terms appear and we are left (formally) with

$$\frac{\mathrm{d}z^k(t)}{\mathrm{d}t} = -i\omega_k z^k + \sum_{p=1}^{d} B_p^k e^{i\theta_p^k}(X_p(t) - X_p(t-\tau)). \tag{117}$$

An integration of (117) gives

$$z^k(t) = z^k(0)e^{-i\omega_k t} + \sum_{p=1}^{d} B_p^k \int_0^t e^{(-i\omega_k(t-r)+i\theta_p^k)}(X_p(r) - X_p(r-\tau))\,\mathrm{d}r, \tag{118}$$

and thus:

$$\xi^k(t) = \mathrm{Re}\,(z^k(0)e^{-i\omega_k t}) + \sum_{p=1}^{d} B_p^k \int_0^t \cos\Big( -\omega_k(t-r) + \theta_p^k\Big)\Big(X_p(r) - X_p(r-\tau)\Big)\,\mathrm{d}r. \tag{119}$$

The boundary terms introduce exogenous forcing from the time series components, themselves. Due to (119) the dominant effect is manifested as narrowband (about the frequency $\omega_k$) modulations of amplitude that superimpose to the pure oscillation $\mathrm{Re}\,(z^k(0)e^{-i\omega_k t})$.

Since pairs of modulated and oscillatory time series (in phase quadrature and) carried out about one dominant frequency may be generated by **Stuart-Landau (SL) oscillators** [59], we propose closures of Eq. (117) built from SL oscillators. By closure we mean that we want to mimic the effects of the term

$$\sum_{p=1}^{d} B_p^k e^{i\theta_p^k}(X_p(t) - X_p(t-\tau)),$$

in Eq. (117), by means of equations involving only the $z^k(t)$-variables plus stochastic ingredients to be determined. In the next section we propose a solution to this problem exploiting both the properties of DAH spectra (such as summarized by Theorem V.1) and ideas borrowed from the theory of **multilayer stochastic models (MSMs)** for inverse modeling [30].

### B. Multilayer Stuart-Landau models

We turn now to the use of DAHMs and DAHCs for deriving inverse stochastic-dynamic models. The purpose is to show that these objects allows for a reduction of the data-driven modeling effort to the learning of elemental MSMs [30] stacked per frequency. These elemental models fall into the class of networks of linearly coupled Stuart-Landau oscillators [32], that may include memory terms as described below; see also [60]. In a certain sense the goal is to show that, given a sequence (of possibly partial) observations issued from a dynamical model, the DAHCs allow for recasting these observations so that they can be modeled within a universal parametric family of simple stochastic models, provided that, as mentioned earlier, the decay of correlations are sufficiently well resolved[61]. We refer to Sect. VIII for examples supporting this statement, and describe hereafter the aforementioned elemental models.

Let $(x_j, y_j)$ be a pair DAHCs such as shown, for instance, in the third column of Fig. 9 or Fig. 12. Except at the zero-frequency and as motivated in Sect. VII A above, DAHCs are typically oscillatory time series that are narrowband about a dominant frequency $f$ and whose oscillations are modulated in time. Such a pair of time series can be reasonably expected to be modeled via nonlinear oscillating systems near a Hopf bifurcation to capture the frequency $f$, and subject to noise, to capture the amplitude modulations; see [59]. We assume that we have a collection of $d$ such pairs, each associated with the same frequency $f$ and exhibiting modulations occurring in different time locations.

If the signal is made of a concatenation of $d$ different time series such as for the L96 flow, such a pair has a simple interpretation: it provides the manifestation, in the time domain, of the frequency $f$ contained in the system's flow as captured by the corresponding pair of DAHMs.

For a narrowband pair of DAHCs, $(x_j(t), y_j(t))$, given at a frequency $f \neq 0$, we envision thus its modeling by some stochastic perturbation of a Stuart-Landau model, namely by

$$\dot{z} = (\mu + i\gamma)z - (1 + i\beta)|z|^2 z + \varepsilon_t, \;\; z \in \mathbb{C}, \tag{120}$$

where $z(t) = x_j(t) + iy_j(t)$, and $\mu, \gamma$ and $\beta$ are real parameters to be estimated. The system is here driven by some "reddish noise," $\varepsilon_t = (\varepsilon_j^x, \varepsilon_j^y)$, to be estimated, also, from the time history of $z(t)$. The issue is that the pairs, when modeled this way, would fail in reproducing the phase coherences across the channels, that would be contained in the signal such as in the case of a wave propagation for instance. The estimation and simulation of $\varepsilon_t$ requires thus a special care due to the presence of typical non-zero cross correlations among the channels.

The MSM framework [30] allows us to address this issue, while keeping the advantages of the modeling framework provided by Eq. (120). It allows, as we will see in applications (see Sect. VIII), for providing stochastic-dynamic inverse models able to respect the phase coherence of the collective behavior of the DAHC pairs $(x_j(t), y_j(t))$ at a given frequency $f \neq 0$. In that respect, one proposes to model each such a DAHC pair by the following type of MSMs, which in the case of two layers used to model the main noise residual, $\varepsilon_t$, can be written as the following system of SDEs:

$$
\begin{aligned}
\dot{x_j} &= \beta_j(f)x_j - \alpha_j(f)y_j + \sigma_j(f)x_j(x_j^2+y_j^2) + \sum_{\substack{i\neq j\\ i\in\mathscr{J}(f)}} b_{ij}^x(f)x_i + \sum_{\substack{i\neq j\\ i\in\mathscr{J}(f)}} a_{ij}^x(f)y_i + \varepsilon_j^x,\\
\dot{y_j} &= \alpha_j(f)x_j + \beta_j(f)y_j + \sigma_j(f)y_j(x_j^2+y_j^2) + \sum_{\substack{i\neq j\\ i\in\mathscr{J}(f)}} a_{ij}^y(f)x_i + \sum_{\substack{i\neq j\\ i\in\mathscr{J}(f)}} b_{ij}^y(f)y_i + \varepsilon_j^y,\\
\dot{\varepsilon}_j^x &= K_{11}^j(f)x_j + K_{12}^j(f)y_j + M_{11}^j(f)\varepsilon_j^x + M_{12}^j(f)\varepsilon_j^y + \zeta_j^x,\\
\dot{\varepsilon}_j^y &= K_{21}^j(f)x_j + K_{22}^j(f)y_j + M_{21}^j(f)\varepsilon_j^x + M_{22}^j(f)\varepsilon_j^y + \zeta_j^y,\\
\dot{\zeta}_j^x &= L_{11}^j(f)x_j + L_{12}^j(f)y_j + P_{11}^j(f)\varepsilon_j^x + P_{12}^j(f)\varepsilon_j^y + N_{11}^j(f)\zeta_j^x + N_{12}^j(f)\zeta_j^y +\\
&\qquad Q_{11}^j(f)\dot{B}_1^j + Q_{12}^j(f)\dot{B}_2^j + \sum_{\substack{i\neq j\\ i\in\mathscr{J}(f)}} \sum_{k=1}^{2} Q_{1k}^i(f)\dot{B}_k^i,\\
\dot{\zeta}_j^y &= L_{21}^j(f)x_j + L_{22}^j(f)y_j + P_{21}^j(f)\varepsilon_j^x + P_{22}^j(f)\varepsilon_j^y + N_{21}^j(f)\zeta_j^x + N_{22}^j(f)\zeta_j^y +\\
&\qquad Q_{21}^j(f)\dot{B}_1^j + Q_{22}^j(f)\dot{B}_2^j + \sum_{\substack{i\neq j\\ i\in\mathscr{J}(f)}} \sum_{k=1}^{2} Q_{2k}^i(f)\dot{B}_k^i.
\end{aligned}
\tag{MSLM}
$$

Such a model is named a **multilayer stochastic Stuart-Landau models** (MSLM). In (*MSLM*), the set $\mathscr{J}(f)$ corresponds of the set of indices $i$ in $\{1,\cdots,dM'\}$ for which the DAHC pair, $(x_i, y_i)$, corresponds to a DAHM pair associated with the frequency $f$. Note that by construction the index $j$ belongs to $\mathscr{J}(f)$, which consists of $d$ elements; see (59) in Theorem V.1.

The $B_k^l$'s with $k$ in $\{1,2\}$ and $l$ in $\mathscr{J}(f)$ form $2d$ independent Brownian motions. At a given frequency $f$, the $d$ pairs are linearly coupled as indicated by the terms in the summation terms appearing in the $x_j$- and $y_j$-equations. As depicted in the schematic of Fig. 6, the main and last layers (here the $x_j$-,$y_j$- and $\zeta_j^{x/y}$-equations) are the only layers in which these pairs are explicitly coupled. Without these coupling terms, we are left with a collection of $d$ uncoupled normal forms of a Hopf bifurcation perturbed by a "reddish" noise, which are the analogue in real domain of the aforementioned Stuart-Landau model (120). The linear coupling at the main level is here aimed to enhance some coherence between the stochastic nonlinear oscillations.

The noise residual $(\varepsilon_j^x, \varepsilon_j^y)$ in the main layer of (*MSLM*) is modeled in the intermediate layer by means of linear dependences involving only $(\varepsilon_j^x, \varepsilon_j^y)$ on one hand, and the $j$-th pair $(x_j, y_j)$, on the other. The noise residual $(\zeta_j^x, \zeta_j^y)$ of that intermediate layer is then modeled still via linear dependences that now include this noise itself, the noise residual of the previous layer, and the DAHC pair $(x_j, y_j)$. The terms in the summations appearing in the $\zeta_j^x$- and $\zeta_j^y$-equations of the last layer involve the other pairs as mentioned earlier. The presence of such terms is aimed to take into account—at the noise level—of cross- correlations between the $d$ DAHC pairs corresponding to a same frequency. In practice, the corresponding coefficient $Q_{lk}^j$ ($l,k \in \{1,2\}$, $j \in \{1,\cdots,d\}$) are estimated from the vector of size $2d$ formed after the concatenation of the last level residuals for each pair.

Note that for the frequency $f = 0$, and as mentioned in Sect. VI D, there are exactly $d$ DAHMs that are not paired, and they are modeled by a linear MSM, following [30]; see also [62, Eq. (9)].

An MSLM is not limited to the two-layer case shown above for the modeling of the main noise residual $(\varepsilon_j^x, \varepsilon_j^y)$. Following [30], extra layers may be added, with intermediate layers that only involve a linear dependence on the $j$-th DAHC pair $(x_j, y_j)$ and on the corresponding noise residuals from the previous layers. In other words, the extra layers must respect the overall structure of (*MSLM*), with coupling between the other pairs (associated with the same frequency $f$) only present in the main and last layers. The addition of an extra layer would thus correspond to the addition of another intermediate plane in the schematic of Fig. 6, with only vertical connections to the rest of the structure.

The procedure of adding layers is controlled by a stopping criterium as described in [30, Appendix A] which allows thus, in principle, for a different number of layers per frequency. In general, the layers with their surrogate variables, here $(\varepsilon_j^x, \varepsilon_j^y)$ and $(\zeta_j^x, \zeta_j^y)$, allow for a learning of the temporal correlations contained in the main noise residual $\varepsilon_t$, as well as of its dependences on the history of $(x_j(t), y_j(t))$ often required for an appropriate modeling of $\varepsilon_t$; see [30] for a general discussion on this topic. Typically, the more layers, the more complex are these correlations and memory dependences; see [30, Proposition 3.3], and [30, Sect. 7], [62] for examples. Sections VIII A 3 and VIII B 3 below, illustrate also this point.

### C. Coupling across the frequencies and practical considerations

As emphasized in (*MSLM*), the coefficients to be learned, depend on the frequency $f$. An MSLM seeks for the modeling of the DAHC pairs $(x_j, y_j)_{j\in\mathscr{J}(f)}$ associated with the frequency $f$, such as extracted from the original (or preprocessed[63]) dataset by application of the DAH decomposition. Once this extraction is computed, an MSLM can be learned in parallel for each frequency of interest, by application of, e.g., simple successive regressions. In practice, linear constraints must be incorporated in these regressions such as $\sigma_j(f) \le 0$ to ensure stability, as well as antisymmetry, to estimate the coefficients $\beta_j(f)$ and $\alpha_j(f)$.

Once all the resulting (few) MSLM coefficients have been estimated, the parallelization advantages of the MSLMs across the frequencies, relies also on its simulation: no need of extra coupling across the frequencies is indeed needed other than running the MSLMs by the same noise realization, for each frequency. The simulation procedure can be thus also parallelized at the convenience of the modeler to obtain the collective behavior of the simulated DAHCs. These simulated DAHCs can be then multiplied with the corresponding DAHMs such as described in Appendix A, to obtain various possible frequency-band modeling of the original dataset.

Of course, one condition to the success of the resulting modeling approach relies on the ability of the DAH decomposition in extracting modulated DAHC pairs that come in nearly phase quadrature, and that are narrowband in the frequency domain, an MSLM being naturally dedicated to the modeling of such time series. The tradeoff between a good resolution of the decay of temporal correlations and the temporal length of the available dataset, seem to play a key role in that respect. Various conducted experiments (not shown) have indeed corroborated the idea that the DAHC pairs are expected to be narrowband and to come in nearly phase quadrature when the decay of temporal correlations has been sufficiently resolved.

In what follows the combination of (i) the extraction of DAHCs by the DAH method of Sect. VI D, and (ii) their modeling according the the MSLM approach described above, is named a **DAH-MSLM modeling**. We present in the next subsections applications to examples for which this criterion is met and a successful DAH-MSLM modeling is obtained in each case.

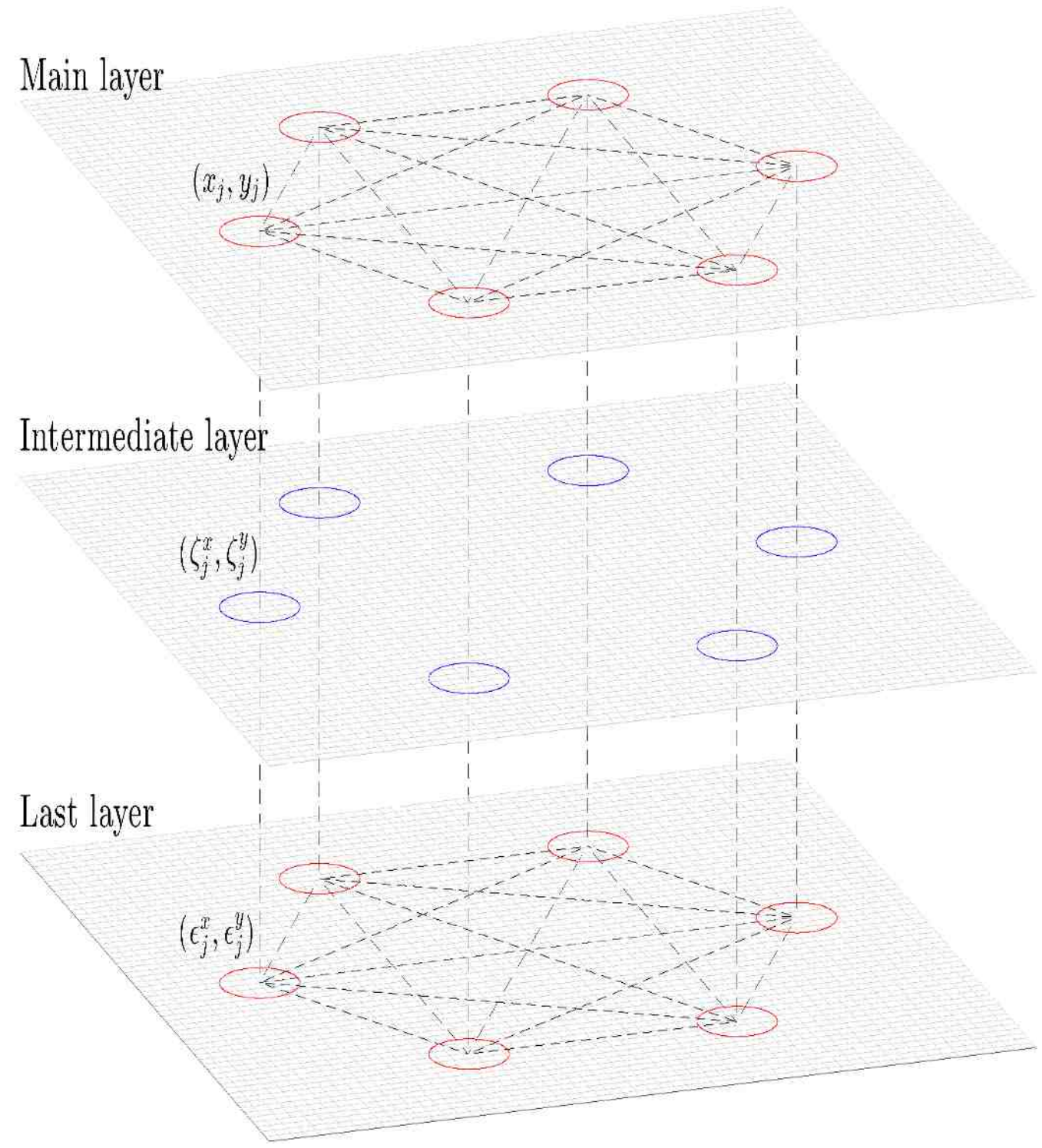

FIG. 6: **Schematic of Eq. (*MSLM*) at a given Fourier frequency $f$.** Each circle represents either a DAHC pair $(x_j, y_j)$ or the corresponding pairs of (successive) noise residuals present in Eq. (*MSLM*). The vertical or horizontal dashed lines connecting two pairs indicate the existence of a mutual dependence between these pairs in the formulation of Eq. (*MSLM*).

## VIII. APPLICATIONS TO INVERSE STOCHASTIC-DYNAMIC MODELING

We apply hereafter the DAH-MSLM inverse modeling approach to two complex flows, obtained respectively from a chaotic ODE system and an SPDE, with linear drift.

### A. Applications to the Lorenz 96 model

#### 1. *The dynamical model*

In this section we apply the DAH approach to the so-called Lorenz 96 (L96) model [64]. The evolution of equations of the L96 model, in its simplified version[65], is written as

$$\frac{dX_k}{dt} = X_{k-1}(X_{k+1} - X_{k-2}) - X_k + F_1, \;\; k = 1, \cdots, K, \quad (121)$$

supplemented with the boundary conditions

$$X_{k-K} = X_{k+K} = X_k. \quad (122)$$

In absence of forcing and dissipation, the sum of the squares of the variables (the energy of the system) is conserved. This model has become popular as a toy model of the atmosphere to test various methods in the weather and climate science as displaying a large set of chaotic behaviors with various mixing properties, as $K$ and $F_1$ is varied.

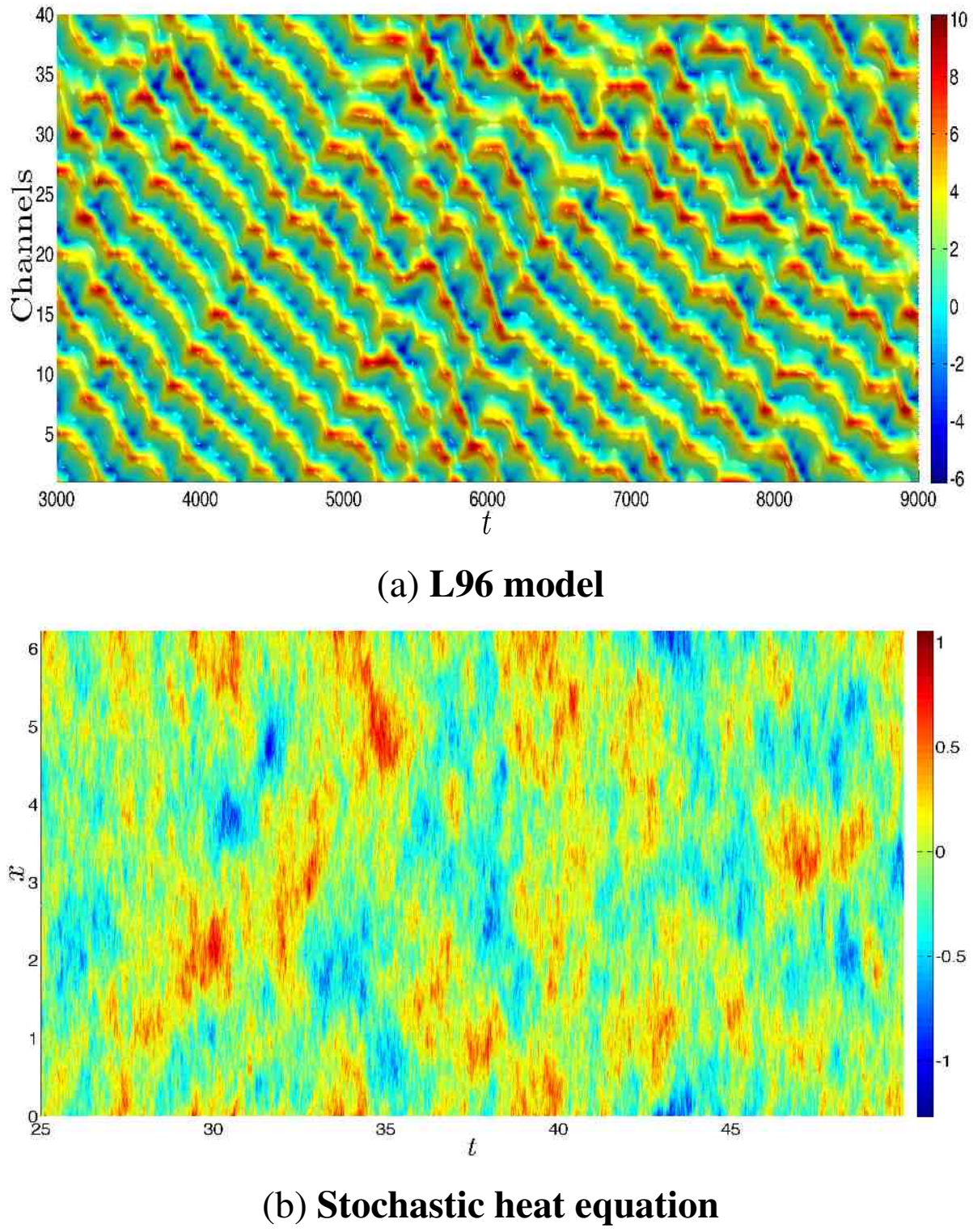

(a) **L96 model**

(b) **Stochastic heat equation**

FIG. 7: **Spatio-temporal fields.** Panel (a) shows a time evolution the Lorenz 96's flow solving Eqns. (121)-(122), and panel (b) shows the solution field of the stochastic heat equation (123). See text for the parameter values and the numerical schemes.

### 2. *DAH spectra, DAHMs and DAHCs*

The model is set to a standard chaotic regime corresponding to $K = 40$ and $F_1 = 6$, and is integrated by using a 4th-order Runge-Kutta scheme with a time step $\delta t = 5 \times 10^{-3}$. Excluding transient behavior and subsampling every 10th time step, we store $N = 11900$ time instances of the $K$-dimensional vector $\mathbf{X}(t) = (X_1(t), \cdots, X_K(t))$.

As shown in Fig. 7-(a), the L96 flow is dominated by irregular spatio-temporal traveling wave-like patterns. We estimated the cross-correlations $\rho^{i,j}(t)$ between the variables $X_i(t)$ and $X_j(t)$ for $i, j$ in $\{1, \cdots, 40\}$, and we formed the $(2M-1) \times (2M-1)$-matrix $\mathfrak{C}$ given in (105), with $M = 160$. The corresponding absolute values of its eigenvalues (ranked per frequency)—forming thus the corresponding DAH power spectrum—are shown in Fig. 8-(a). The corresponding DAH phase spectrum is shown in Fig. 8-(b).

Interestingly, these spectra as well the underlying DAHMs share some common features with those displayed in the case of a pure traveling wave shown in Sectns. VI E 1 and VI E 2. First, the DAHMs shown in the 1st and 2nd columns of Fig. 9 share clearly common patterns with those shown in the respective columns of Fig. 5.

Second, a phase synchronization between the channels of a DAHM as $f$ approaches the Nyquist frequency, is also observed in Fig. 8-(b) as in Fig. 3-(b).

Third, this synchronization is in contrast with the disordered distribution of phases shown in a frequency band $0 < f < f_c$ with $f_c \approx 0.1$; frequency band for which the most energetic pair of DAH eigenvalues (for the quadratic form $\mathscr{Q}_\delta$) are located; see Fig. 8-(a).

Nevertheless, the DAH power spectrum exhibits much broader peaks in the case of the L96 flow than in the case of a pure traveling wave; compare Fig. 8-(a) with Fig. 3-(a).

Finally, the most noticeable difference is shown by the DAHCs. Whereas the latter were all periodic in the case of traveling wave as expected (see (115)), they are—as predicted by (119)[66]—clearly modulated in amplitude and narrow-band in frequency in the case of the L96 flow; compare the 3rd column of Fig. 5 with that of Fig. 9. Furthermore the pairs of DAHCs are typically constituted by time series (for $f \neq 0$) that are nearly in phase quadrature; see the panels corresponding to $f = 0.016$, $f = 0.031$ and $f = 0.047$ in the 3rd column of Fig. 12. We are thus in a situation described in Sect. VII B. The next section addresses the inverse-modeling problem of the L96 flow, following the approach described in Sect. VII B.

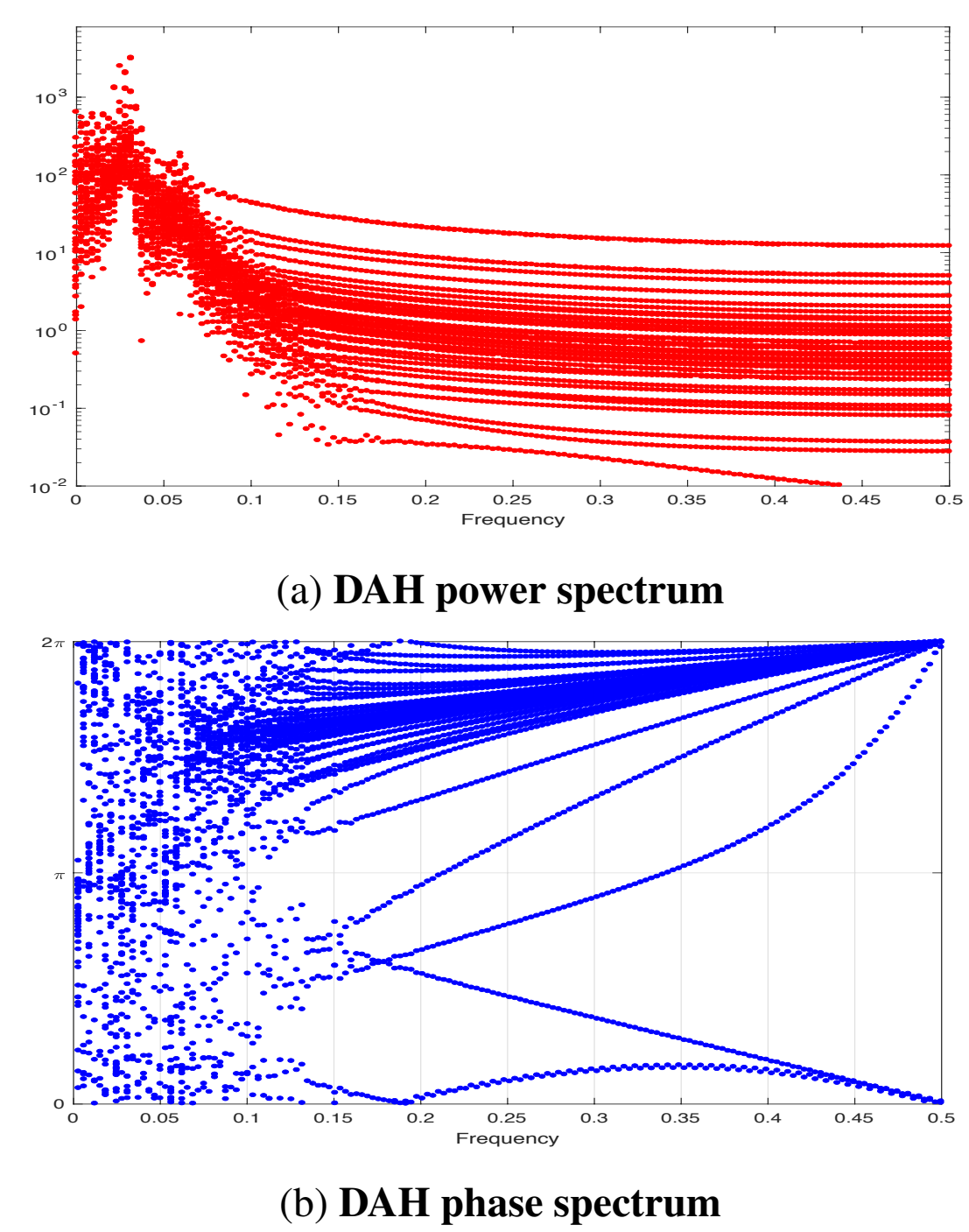

(a) **DAH power spectrum**

(b) **DAH phase spectrum**

FIG. 8: **Multidimensional DAH power and phase spectra for the L96 flow.** We chose $M = 160$ to estimate the cross correlations and form the corresponding the matrix $\mathfrak{C}$ given in (105). The DAH power spectrum and DAH phase spectrum are obtained as described in Sect. VI D. Each non-zero frequency is associated with 40 eigenpairs that correspond to the total of $d = K = 40$ variables in the L96 model.

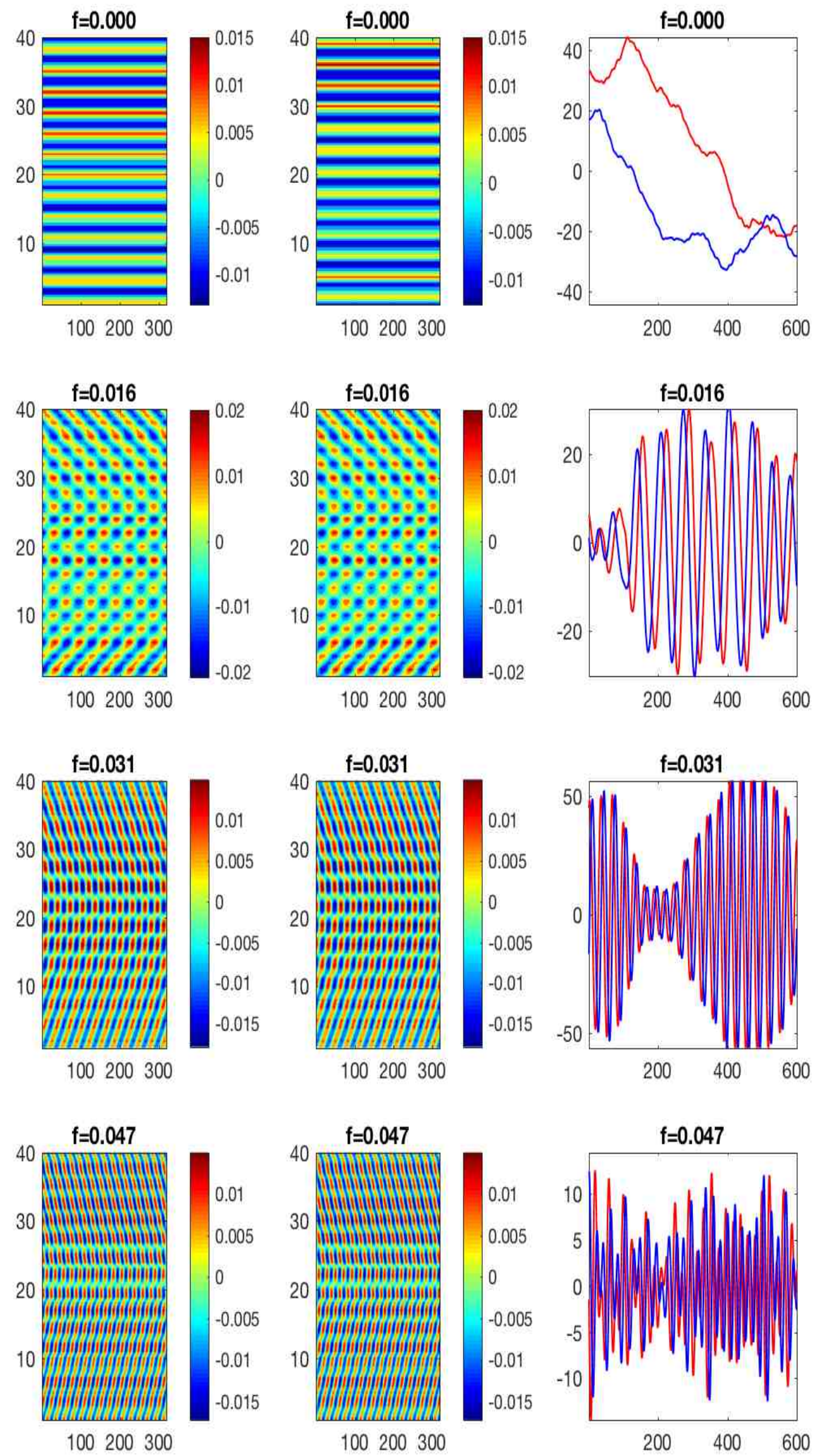

FIG. 9: **DAHMs and DAHCs from the L96 flow.** Same as in Fig. 4 but for the L96 flow shown in Fig. 7-(a). Here $M' = 2M - 1$ (here $M = 160$), and $d = 40$. The DAHMs and DAHCs shown here correspond, for each indicated frequency $f$, to the dominant pair of DAH eigenvalues located above $f$ in Fig. 8-(a).

### 3. *DAH-MSLM modeling of the Lorenz 96 model*

In this section we report on the inverse modeling skills obtained by the DAH-MSLM approach of Sect. VII. For that purpose, the original field obtained by integrating Eq. (121) as mentioned above, has been sampled every 10 time steps in time. The frequency band that is modeled is from $f = 0$ up to $f = 6.29 \times 10^{-2}$, which represents a fraction of the variance of 86%, and 20 Fourier frequencies according to the resolution set by $M' = 2 \times 160 - 1 = 319$.

The corresponding DAHCs are modeled for each of these frequencies according to an MSLM of the form described in VII B, in which the number of layers per frequency is 11, including the main one. After multiplication by the the corresponding DAHMs, the simulated DAHCs allow us, as described in Appendix A, to form a **simulated DAH-MSLM field**. The resulting field is shown in Fig. 10 and exhibits striking common features with the original field, shown on the same figure for a better comparison.

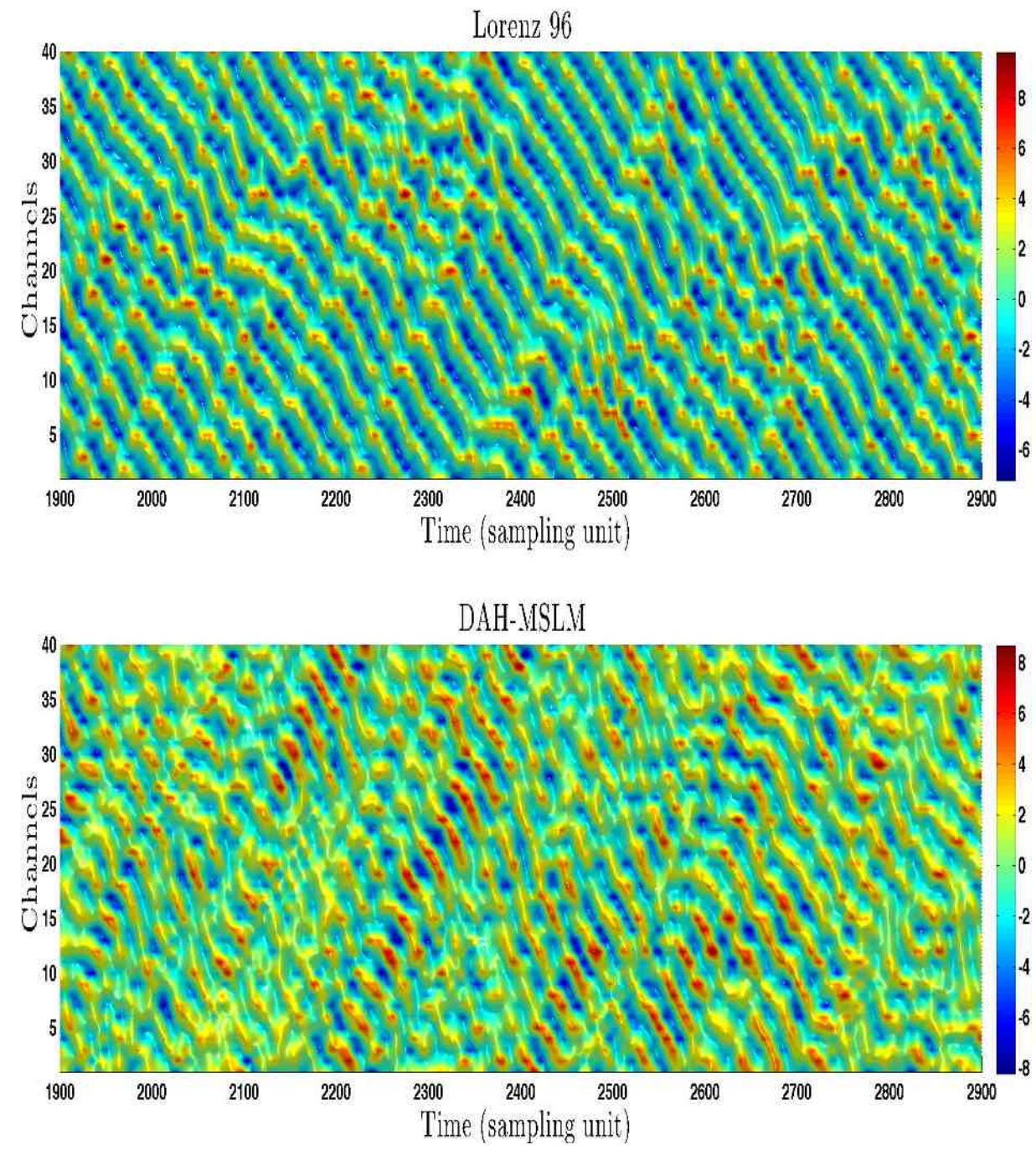

FIG. 10: **The L96 flow and its simulated DAH-MSLM field.** The L96 flow shown in the top panel corresponds to the sampling units used at the learning stage, i.e. every 10th time steps; see text for more details.

Interested by the ability of the simulated DAH-MSLM field to reproduce the time-variability of the L96 flow, we computed for each of these fields the autocorrelation function (ACF) (in time) of each channel, followed by an averaging of the resulting individual ACFs over the channels. The result is shown in Fig. 14-(a) and as one can see the ability of the simulated DAH-MSLM field in reproducing the time-variability of the L96 flow is very satisfactory, according to this metric.

It should be emphasized that beyond the embedding window size, here $M' = 319$, these skills deteriorate due in particular to the periodization involved in the determination of the matrix $\mathfrak{C}$, as in any other method exploiting time-embedding techniques.

So far we have been dealing in applications with deterministic systems. The next section illustrates that striking inverse modeling skills can also be achieved by applying the DAH-MSLM approach for the inverse modeling of stochastic systems, still with the same elemental MSLMs as described in Sect. VII B.

### B. Applications to a stochastic heat equation

#### 1. The dynamical model

The dataset analyzed hereafter is produced from a simple stochastic partial differential equation (SPDE), namely the following stochastic heat equation posed on the interval $[0,L]$ ($L>0$) with periodic boundary conditions:

$$\begin{aligned} du(t,x) = (-\lambda u(t,x) + D\partial_{xx}^2 u(t,x))dt + \sqrt{2\lambda}\, dW(t,x), \\ x \in \mathbb{T}_L = \mathbb{R}/(L\mathbb{Z}). \end{aligned} \tag{123}$$

Here $\lambda$ and $D$ are positive parameters, and $dW$ denotes a "space-time" white noise process on the circle $\mathbb{T}_L$, roughly speaking a Gaussian random field with spatio-temporal correlations for $t,s \geq 0$, given by

$$\mathbb{E}(\, dW(t,x)\, dW(s,y)) = \delta_0(t-s)\delta_0(y-x), \;\; x,y \in \mathbb{T}_L.$$

There are (at least) two approaches to frame mathematically the notion of space-time white noise and related concepts of solutions; see e.g. [67] for a comparison. Both, the random-field approach based on Green functions [68] and the functional analysis approach [69] [70], give the same notion of weak/mild solutions, and stochastic heat equations such as Eq. (123) serve as a fundamental object to study various properties or formulas about such solutions; see e.g. [71–74].

The functional analysis approach teaches us that the stochastic process $u(t,\cdot)$ solving Eq. (123) is actually an **infinite-dimensional Ornstein-Uhlenbeck process (OU)** in the Hilbert space $\mathscr{H} = L^2(\mathbb{T}_L)$; see [71]. Our goal is to analyze the DAH power and phase spectra of the corresponding (random) spatio-temporal field, $u(t,x)$, and see whether these spectra provide an intuitive characterization of such an OU process.

For that purpose, we first simulate Eq. (123) as follows. Over a grid of mesh size $\delta x = L/N_x$, the discrete approximation $u_j^n$ of $u(n\delta t, j\delta x)$ is obtained by finite differences for which the noise term is approximated at each time step $n\delta t$ (via an explicit Euler-Maruyama scheme) by $\xi_n\sqrt{\delta t}$ where $(\xi_j^n)_{j=1}^{N_x}$ denotes an $N_x$-dimensional vector of random variables each drawn independently with respect to $1 \leq j \leq N_x$ and $n$ (except $\xi_1^n = \xi_{N_x}^n$) from the standard normal distribution $\mathscr{N}(0,1)$.

We chose $N_x = 2^7$, $L = 2\pi$, $D = 2\times 10^{-1}$ and $\delta t = 10^{-3}$ for our simulations. The initial data is taken to be

$$u_0(x) = 10^{-1}\sqrt{\frac{2}{L}}\cos\left(\frac{2\pi x}{L}\right). \tag{124}$$

Without any surprise, the integration reveals striking differences in term of patterns exhibited by the field $u(t,x)$ compared to those obtained from the nonlinear example of Section VIII A; compare panel (b) with panel (a) in Fig. 7.

#### 2. DAH spectra, DAHMs and DAHCs

Excluding transient behavior and subsampling every $8\delta t$, we store $N = 37500$ time instances of the $d$-dimensional vector $\mathbf{X}(t) = (X_1(t),\cdots,X_d(t))$, in which the $X_p$'s are also obtained by subsampling every $4\delta x$ the original field $(u_j^n)$, i.e. here $d = 32$. The resulting coarse-grained SPDE field is shown in the top panel of Fig. 13. We built up the matrix $\mathfrak{C}$ (105) whose Hankel blocks are given in (104) with $M = 80$, and $1 \leq i,j \leq d$. Diagonalizing $\mathfrak{C}$ gives the multidimensional DAH spectra as explained in Sect. VI D. The corresponding DAH power (resp. DAH phase) spectrum for the coarse-grained SPDE is shown in Fig. 11-(a) (resp. Fig. 11-(b)).

As it can be observed by comparison with Fig. 8-(a), the DAH power spectrum shown in Fig. 11-(a) does not exhibit any broadband peak, and is actually composed by curves of DAH eigenvalues that each exhibits a decay from low to high frequencies, reminiscent with that of a scalar OU process. What is remarkable is that the DAH power spectrum provides here a multidimensional picture of what one should expect for an OU process obtained from a standard 1D Langevin equation, showing thus here again the relevance of the quadratic form $\mathscr{Q}_\delta$, as an energy characterizing the distribution of eigenvalues of $\mathfrak{C}$ at a given frequency; see again the end of Sect. VI D and Remark V.1-(ii).

As another noticeable feature, the phase synchronization such as observed on the DAH phase spectra of the L96 flow and the traveling wave examples, is no longer present for this stochastic example, and is instead replaced in Fig. 11-(b) by the presence of some (vaguely visible) bands within a diffuse background. Going back to the interpretation of the phase spectrum discussed at the end of Sect. V in relation with (81), the presence of such a diffuse background, even for $f$ approaching the Nyquist frequency $f_N$, can be attributed to the presence of non-vanishing cross-spectral terms $\widehat{\rho^{p,q}}(f)$ leading to a non-negligible $R_p(f)$ given in (79). This feature is consistent with the intuitive idea that the noise-term in Eq. (123) is white and thus forces evenly the resolved frequency band. It should be emphasized that such a feature has been also observed in other conducted experiments on datasets issued from various stochastic systems driven by a white and additive noise.

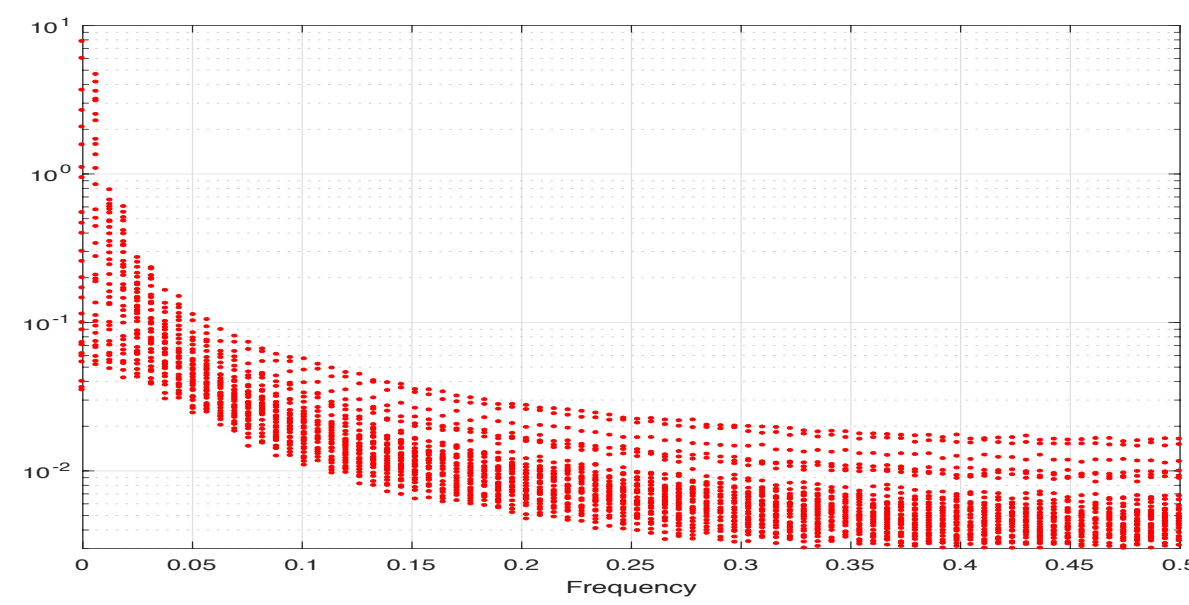

(a) **DAH power spectrum**

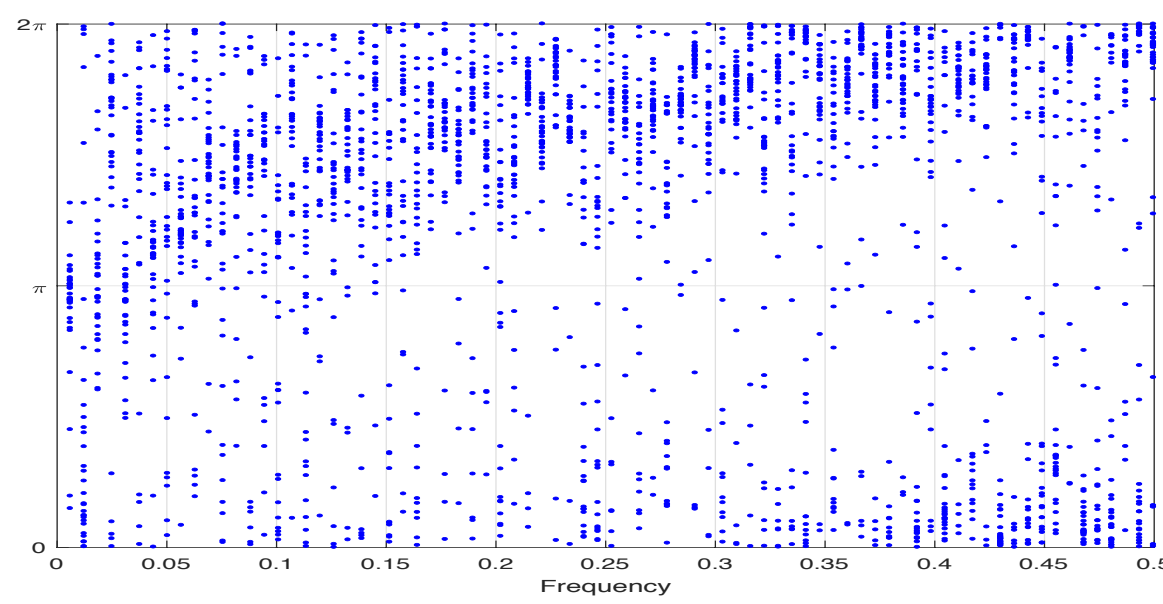

(b) **DAH phase spectrum**

FIG. 11: **Multidimensional DAH power and phase spectra of the SPDE flow.** We chose $M = 80$ to estimate the cross correlations and form the corresponding matrix $\mathfrak{C}$ given in (105). The DAH power spectrum and DAH phase spectrum are obtained as described in Sect. VI D. Each non-zero frequency is associated with 32 eigenpairs that correspond to the total of $d = 32$ spatial channels obtained by coarse-graining the original dataset as described in text.

Finally, each pair of DAHCs associated with a frequency $f$ (with $f \neq 0$) is—as for the nonlinear example of Sect. VIII A—constituted by time series modulated in amplitude, narrowband about the frequency $f$, that are nearly in phase quadrature; see the panels corresponding to $f = 0.013$, $f = 0.025$ and $f = 0.038$ in the 3rd column of Fig. 12. We are thus here again in a favorable situation to apply to these DAHCs the MSLM modeling framework of Sect. VII B. The next section examine the corresponding modeling skills.

### 3. *Coarse-grained DAH-MSLM modeling*

Recall that the original SPDE field has been sampled every 8-th time step, and every 4-th spatial increment $\delta x$ in space; the resulting coarse-grained field is shown on the top panel of Fig. 13. The frequency band that is modeled is from $f = 0$ up to $f = 6.33 \times 10^{-2}$, which represents a fraction of the variance of 86%, and corresponds to 10 first Fourier frequencies according to the resolution set by $M' = 2 \times 80 - 1 = 159$.

The corresponding DAHCs are here again modeled for each of these frequencies according to an MSLM of the form described in Sect. VII B, in which the number of layers per frequency is 4, including the main one. The resulting simulated DAH-MSLM field is shown in Fig. 13 and exhibits striking common features with the coarse-grained SPDE field shown on the same figure.

At a statistical level, the simulated DAH-MSLM field (red curve in Fig. 15) captures remarkably well the distribution of energy across the spatial scales when compared with the coarse-grained version of the SPDE field (blue curve in Fig. 15) from which the parameters the DAH-MSLM have been estimated. Note that the blue curve in Fig. 15 shows this energy distribution for the coarse-grained SPDE field as projected onto the first pairs of DAHMs corresponding to the first 10 Fourier frequencies. The black curve shows in the same figure, the energy distribution for the full coarse-grained SPDE field.

Finally, as in the case of the L96 flow, in order to assess the ability of the simulated DAH-MSLM field to reproduce the time-variability contained in the coarse-grained SPDE field, we computed for each of these fields the ACFs of each channel, followed by an averaging of the resulting individual ACFs over the channels. The result is shown in Fig. 14-(b) and as one can see the ability of the simulated DAH-MSLM field in reproducing the (average) decay of correlations is very satisfactory. Noteworthy is the nearly exponential decay of the later, as intuitively expected from an infinite-dimensional OU process.

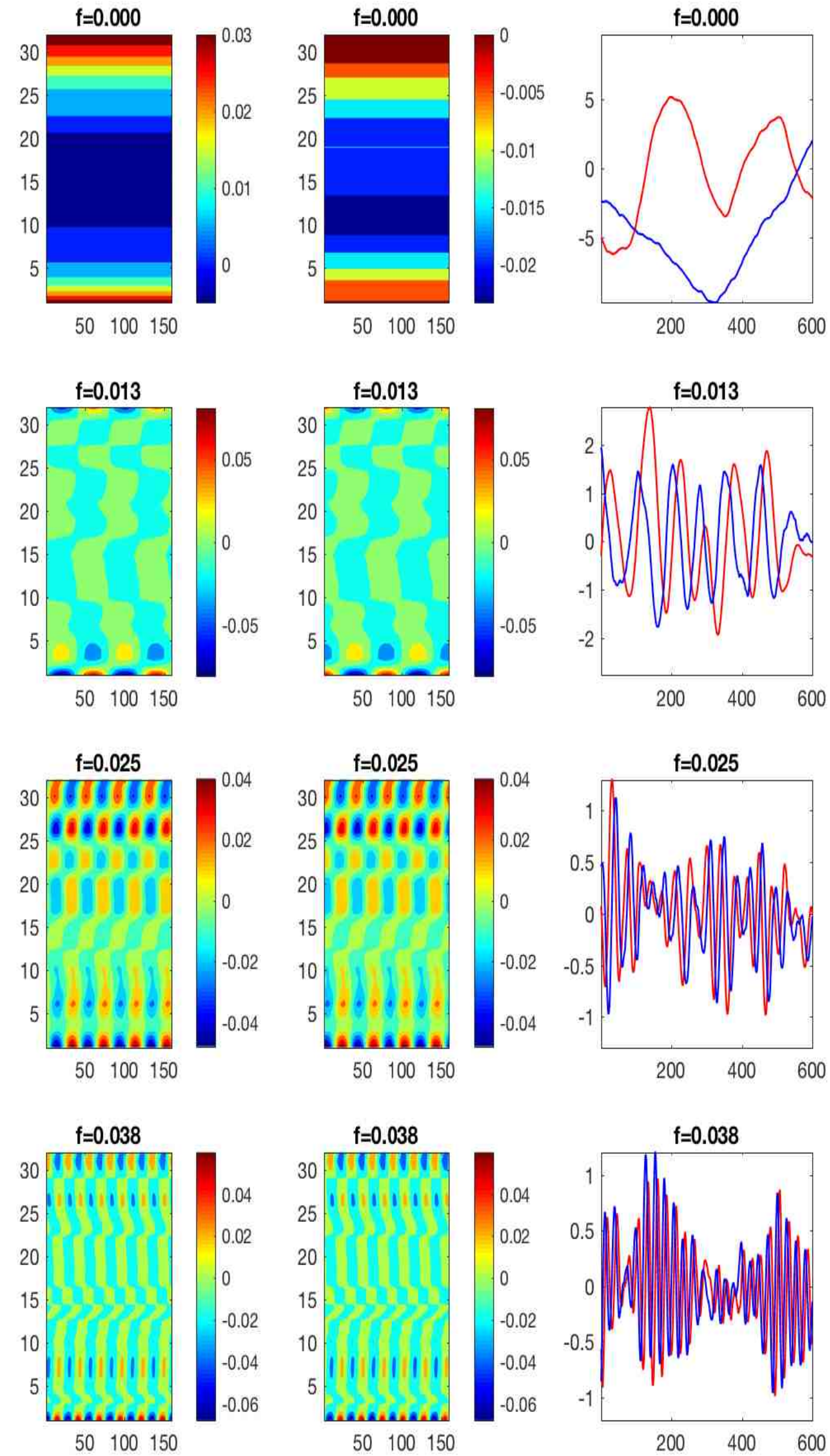

FIG. 12: **DAHMs and DAHCs from the SPDE flow.** Same as in Fig. 4 but for the SPDE flow shown in Fig. 7-(b). Here $M' = 2M - 1$ (here $M = 80$), and $d = 32$. The DAHMs and DAHCs shown here correspond, for each indicated frequency $f$, to the dominant pair of DAH eigenvalues located above $f$ in Fig. 11-(a).

To summarize, the DAH-MSLM inverse modeling framework provides here again efficient emulators of the observed dynamics for many aspects. The most striking feature of the approach is that either for the linear stochastic example examined here or the nonlinear chaotic one of Sect. VIII A, the predictor functions stay the same, namely the monomials used for the class of MSLMs described in Sect. VII B. The data-adaptive character of the DAHMs manifested by the capture of the phase information contained[75] in the multivariate signal to be modeled, constitutes a key element for a successful modeling. It allows indeed for breaking the signal into elementary bricks, here the DAHCs, that fall naturally within the MSLM class of models. In other words, the DAHMs provide a change of basis which transforms systematically the original dataset into time series generated by a same class of SDEs, provided the corresponding DAHCs are modulated in amplitude and narrowband (for $f \neq 0$). Of course many questions arise, as it did for many novel mathematical concepts and tools, when first introduced. The next and final section formulate some of these.

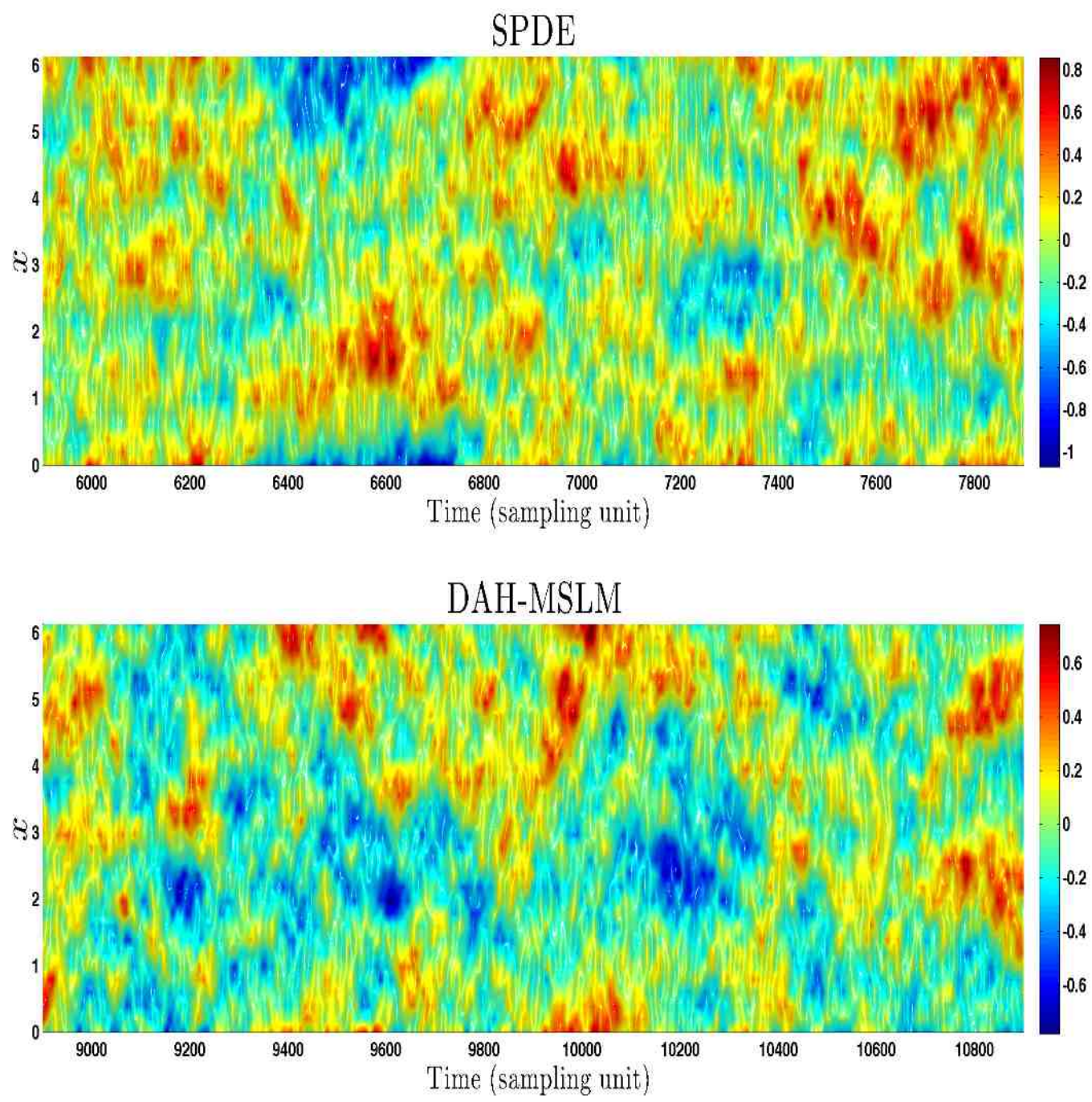

FIG. 13: **The stochastic heat flow and its simulated DAH-MSLM field.** The coarse-grained SPDE field shown in the top panel corresponds to that obtained from sampling units used for the learning stage, i.e. every 8th time steps and 4th spatial mesh size, and after projection onto the pairs of DAHMs corresponding to the first 10 Fourier frequencies; see text for more details.

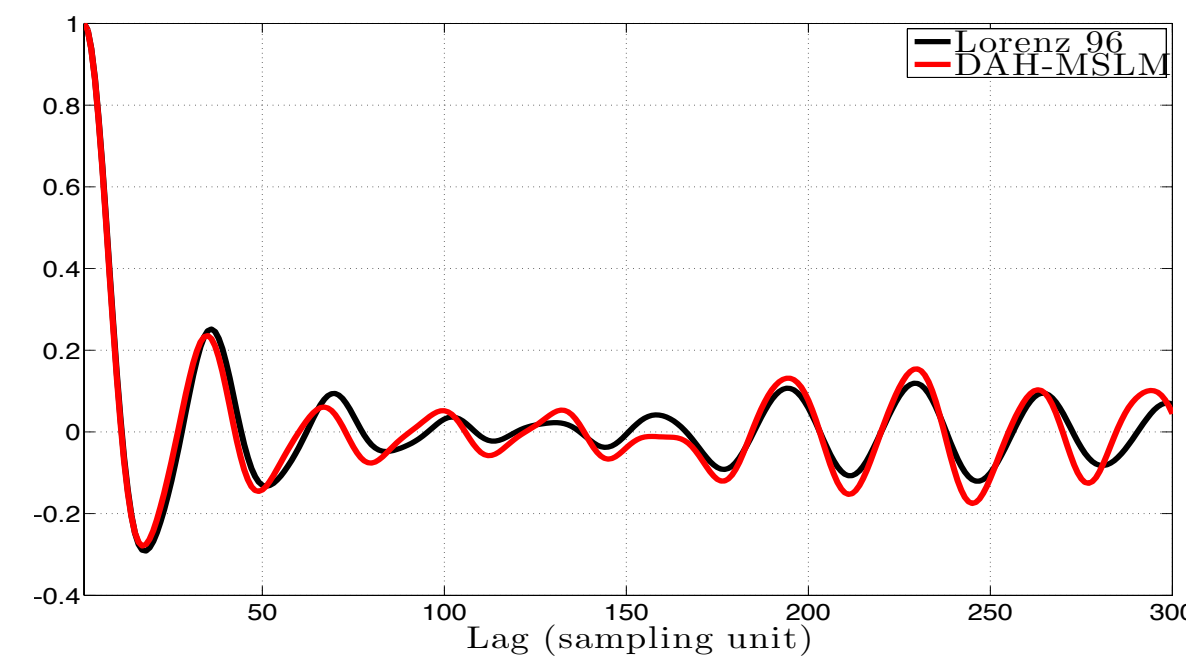

(a) **DAH-MSLM vs Lorenz 96**

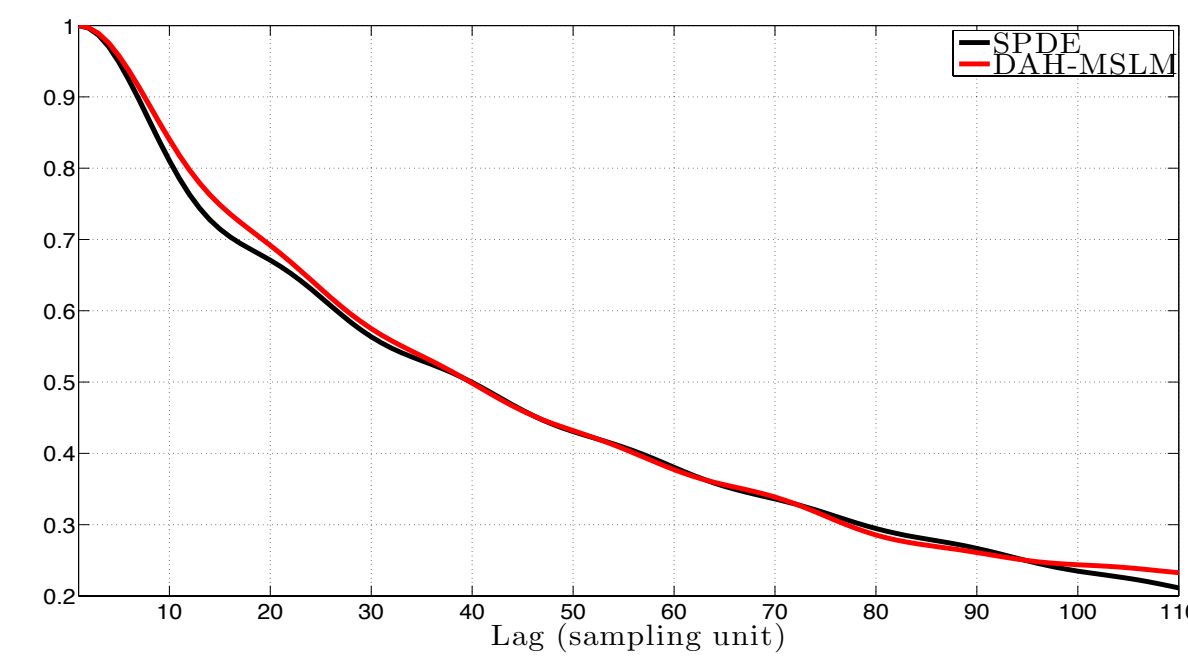

(b) **DAH-MSLM vs SPDE**

FIG. 14: **Averaged autocorrelation functions: modeling skills.** In each panel the black curve is obtained by computing the ACFs of each channel of the corresponding dataset, followed by an averaging of the resulting individual ACFs over the channels. The red curve is obtained by averaging furthermore across the realizations of a given DAH-MSLM simulation.

## IX. CONCLUDING REMARKS

Thus, by means of rigorous harmonic analysis tools, this article provides a natural framework for the representation of complex time-evolving fields by simple time series, the DAHCs, suitable for their modeling within a universal class of models, namely the MSLMs. The framework opens up several possible directions for future research. We outline some of these directions below and comment on recent achievements obtained on non-synthetic datasets.

1. The modeling approach described in this article may benefit from several natural developments. First, terms of the form $\gamma_j(f) y_j(x_j^2 + y_j^2)$ and $\gamma_j(f) x_j(x_j^2 + y_j^2)$ in the $x_j$- and $y_j$-equations, respectively, could be incorporated within the formulation of an MSLM. The introduction of such nonlinear "twist" terms is in part motivated by the analysis conducted in [59]. There, it has been indeed shown that such terms are responsible of certain interactions between the noise and the nonlinear terms that may be important to resolve. The inclusion of such twist terms may be thus relevant for the proper modeling of a DAHC pair, for certain applications.

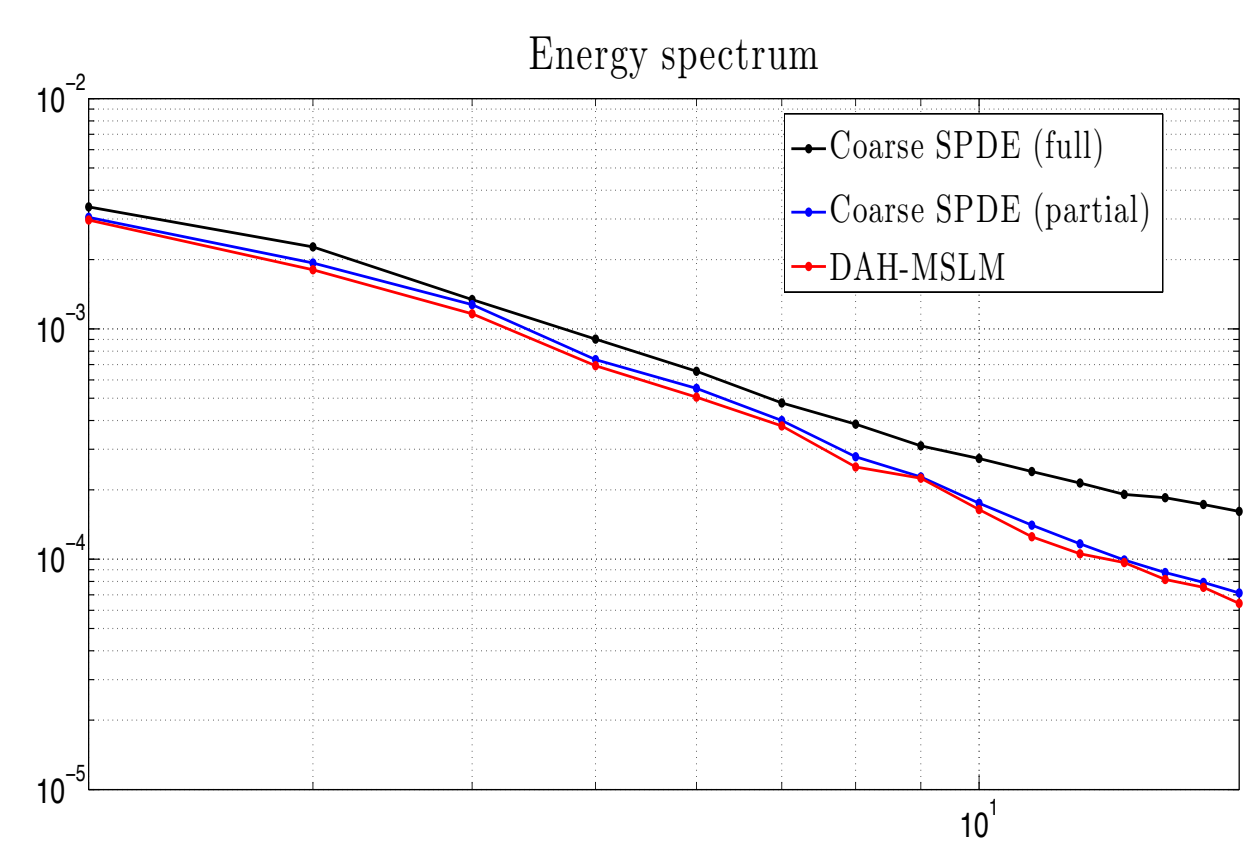

FIG. 15: **Energy spectrum.** Distribution of energy across the spatial scales. The blue curve shows this energy distribution for the coarse-grained SPDE field as projected onto the first pairs of DAHMs corresponding to the first 10 Fourier frequencies; see text for more details. The black curve shows this energy distribution for the full coarse-grained SPDE field (i.e. without projection).

Also, the driving noise of an MSLM may not be required to belong to the class of correlated white noise processes[76] such as in Eq. (*MSLM*). For instance, an MSLMs driven by Lévy $\alpha$-stable processes could be relevant for certain applications such as the modeling of precipitations in climate [77, 78] and other climate fields [79].

2. The use of MSLMs may not be limited to the modeling of pairs of DAHCs as extracted by DAH. For instance, other multivariate time series decomposition techniques such as MSSA [80] extract analogue time series known as ST-PCs, which when coming in pairs that are nearly in phase quadrature, could be also efficiently modeled by an MSLM provided the corresponding pairs are modulated in amplitude while narrowband. We mention that a similar approach was adopted in [62] for the stochastic modeling of decadal variability in ocean gyres. There however, the nonlinear coupling within a frequency bin was absent while linear couplings across the frequencies (falling within a targeted frequency band) were used to model the ST-PCs from the dataset considered in [62].

3. Similarly, coupling across frequencies that go beyond the driving by the same noise realization may also be useful to build up MSLMs, in certain cases. The inclusion into an MSLM of the contribution of another frequency $g$ than a frequency $f$, may be indeed easily achieved by passing the corresponding terms associated with $g$ into (*MSLM*). Of course, as far as the pairs of DAHCs are sufficiently narrowband, such a coupling may be superfluous and in such a case, a "same noise realization coupling"—cheaper in terms of coefficients to be estimated[81]—may be amply sufficient to reach satisfactory modeling skills as shown for the time-evolving fields of Sect. VIII.

4. Unlike the synthetic datasets of Sect. VIII, the DAH-MSLM approach has been applied to various geophysical multivariate datasets exhibiting complex patterns, from "real-world applications." In that respect, the DAH-MSLM modeling of Arctic sea ice concentration has been successfully addressed in [82], providing in particular stable inverse models for much longer time intervals than their training period (1979-2014), beyond the end of the 21st century. Interestingly, these models exhibit an interdecadal variability consistent with proxy historical records, while reproducing accurately the (nonlinear) trend contained in the observations; see [82, Fig. 13].

Prediction experiments using the DAH-MSLM framework have been also recently conducted in [83] for the Arctic sea ice extent (SIE). Here the forecasts, as it would be the case for any other dataset, are made for the DAHCs using the best MSLMs learned during the training period. As a result, retrospective forecasts show that the resulting DAH-MSLM modeling is quite skillful in predicting the September SIE; i.e. for the most challenging SIE value to forecast. Noteworthy are the real-time forecasts conducted in 2016 which have shown very competitive skills as documented here: `https://www.arcus.org/sipn/sea-ice-outlook/2016/post-season`; see [83] for more details.

Finally we mention that the DAH-MSLM framework has been shown to be useful for the data-driven study of the coupling between the solar wind and the magnetosphere as discussed in [84]. This example shows in particular the relevance of the DAH-MSLM framework to provide inverse models of externally forced systems, another important theme for geophysical applications [85, 86].

## ACKNOWLEDGMENTS

The authors want to warmly thank Georg Gottwald for the inspiring discussions at the start of this project. The authors are grateful to Henk Dijstra, David Neelin, and Alexis Tantet for fruitful and inspiring collaborations about the stochastic Hopf bifurcation. Constructive comments by Pavel Berloff, Sergey Kravtsov, Honghu Liu and Eric Simonnet in early stages of this project, are also most appreciated. The authors are also grateful to James McWilliams for several inspiring discussions and constructive comments about this work. Finally, Michael Ghil has the full gratitude of the authors for the numerous stimulating discussions over the years as well as for his unflinching support.

This research was supported by grant N00014-16-1-2073 from the Multidisciplinary University Research Initiative (MURI) of the Office of Naval Research, and by the National Science Foundation grants OCE-1243175 and DMS-1616981.

## Appendix A: DAH coefficients (DAHCs) and frequency band reconstruction

Let us consider a $d$-dimensional time series $\mathbf{X}(t) = (X_1(t), \cdots, X_d(t))$ from which the cross-correlation coefficients $\rho^{i,j}$ between the $i$-th and the $j$-th channels, have been estimated. Once the DAHMs have been determined as eigenvectors of the matrix $\mathfrak{C}$ defined in (105) (see Section VI D), one can use them to determine extract from $\mathbf{X}(t)$ the part of the signal that corresponds to a certain frequency band $[f_s, f_e]$. We describe hereafter how to proceed to do so.

First, note that one can always project the multidimensional time series $\mathbf{X}(t)$ onto the orthogonal set formed by the DAHMs, $W_j$, given as

$$W_j = (E_1^j, \cdots, E_d^j), \;\; j \in \{1, \cdots, dM'\},$$

where $E_p^j$ denotes the $M'$-dimensional discrete approximation of its continuous counterpart given in (61). Sometimes $E_p^j$ is called a **DAHM snippet** or simply a snippet.

After projection onto a DAHM, $W_j$, one obtains thus the following **DAH coefficients (DAHCs)**:

$$\xi_j(t) = \sum_{s=1}^{M'} \sum_{p=1}^{d} X_p(t+s-1) E_p^j(s), \qquad 1 \le s \le M', \; 1 \le t \le N' = N - M' + 1. \tag{A1}$$

Although the DAHCs are not formally orthogonal in time, our numerical experiments indicate that they also exhibit a certain phase-quadrature relationship that seems however to be subordinated to the decay of correlations: roughly speaking, the best the latter are resolved, the more apparent is the phase-quadrature (not shown); see also Sect. VII C.

A partial or full (i.e. without loss of variance) reconstruction of the original data can be then obtained from the DAHCs. Such a reconstruction is obtained from a finite convolution between either a subset or the whole set of DAHCs, and the corresponding $E_p^k$'s.

More exactly, given a $M'$-dimensional DAHM snippet, $E_p^j(s)$ $(1 \leq s \leq M')$, and its corresponding DAHC, $t \mapsto \xi_j(t)$, one can form the **reconstructed component (RC)**, $R_p^j(t)$, given at time $t$ by

$$R_p^j(t) = \frac{1}{M_t} \sum_{s=L_t}^{U_t} A_j(t-s+1) E_p^j(s),\ 1 \leq s \leq M'. \tag{A2}$$

The values of the normalization factor $M_t$, as well as of the lower and upper bound of summation $L_t$ and $U_t$, differ between the central part of the time series and its end points, and are given here as for M-SSA; see [80, Eq. (12)]. Forming the corresponding reconstructed multivariate time series $\mathbf{R}^j(t) = (R_1^j(t), \cdots, R_d^j(t))$, and summing over all of the (resolved) $j$'s allows us to recover the original multivariate time series $\mathbf{X}(t)$. A summation over a subset of $\{1, \cdots, dM'\}$ gives a partial reconstruction of $\mathbf{X}(t)$.

Unlike for M-SSA [80] or PCA, such a partial reconstruction is not necessarily guided by a certain targeted percentage of the variance of $\mathbf{X}(t)$ to be captured, but can also be performed over a targeted frequency band. Due to its harmonic flavor, the DAH approach allows us indeed to consider a notion of **harmonic reconstructed component (HRC)** which consists—for a given channel $p$—of summing the RCs corresponding to a given frequency $f_k$, namely

$$R_p^{f_k}(t) = \sum_{j \in \mathscr{J}(f_k)} R_p^j(t), \tag{A3}$$

where $\mathscr{J}(f_k)$ corresponds to the indices $j$ in $\{1, \cdots, dM'\}$ for which the DAHM snippet, $E_p^j(s)$, is associated with the frequency $f_k$. We know indeed from Theorem V.1 that there is $2d$ eigenvalues corresponding to a frequency $f_k$ and as much DAHMs.

In applications, the DAH approach offers thus a way to determine the contribution within a component $X_p(t)$ (of the multivariate signal $\mathbf{X}(t)$) of a particular frequency band $I_f := [f_s, f_e]$ by simply adding up the HRCs over the frequencies lying in $I_f$.